\documentclass[12pt]{amsart}
\makeatletter
\@namedef{subjclassname@2020}{%
  \textup{2020} Mathematics Subject Classification}
\makeatother
\usepackage{fullpage}
\usepackage{amsthm}
\usepackage{amsmath}
\usepackage{amssymb}
\usepackage{graphicx}
\usepackage{placeins}
\usepackage[bookmarks=true,unicode=true,colorlinks=true,linkcolor=blue,citecolor=blue,urlcolor=blue]{hyperref}
\hypersetup{pdftitle={Behavior of the Caputo time fractional equations as the order alpha tends to zero: the resolvent limit and its stochastic interpretation},pdfauthor={Cheng-Gang Li}}
\newtheorem{thm}{Theorem}[section]
\newtheorem{cor}[thm]{Corollary}
\newtheorem{lem}[thm]{Lemma}
\newtheorem{prop}[thm]{Proposition}
\theoremstyle{definition}
\newtheorem{defn}[thm]{Definition}
\newtheorem{example}[thm]{Example}
\theoremstyle{remark}
\newtheorem{rem}[thm]{Remark}
\numberwithin{equation}{section}
\begin{document}
\title[Caputo time fractional equations as $\alpha\rightarrow 0^+$]{Behavior of the Caputo time fractional equations as the order $\alpha\rightarrow 0^+$: the resolvent limit and its stochastic interpretation}

\author{Cheng-Gang Li}
\address{School of Mathematics, Southwest Jiaotong University, Chengdu 611756, P.R. China}
\email{lichenggang@home.swjtu.edu.cn}

\subjclass[2020]{26A33, 34A08, 35R11, 47D06, 60G51}
\keywords{Caputo fractional derivative, Mittag-Leffler function, resolvent family, inverse stable subordinator, probabilistic potential theory}

\begin{abstract}
We study the behavior, as the fractional order $\alpha\rightarrow0^+$, of the Caputo time derivative and of the solutions of $D_t^\alpha u=Au$ in a Banach space. The solutions converge for every $t>0$ to the time-independent resolvent value $(1-A)^{-1}\phi$. A single bounded fractional resolvent family already forces the sectoriality of $A$ and yields operator-norm convergence with rate $O(\alpha)$ and a full expansion in powers of $\alpha$. The engine is an expansion of the Mittag-Leffler function in its order, $E_\alpha(z)=(1-z)^{-1}+\sum_{k=1}^{K}c_k\alpha^k\mathrm{Li}_{-k}(z)+O(\alpha^{K+1})$ on sectors away from $[1,\infty)$, the corresponding qualitative limit going back to Mittag-Leffler's summability of the geometric series. We also develop stochastic interpretations using inverse $\alpha$-stable subordinators, and establish connections with probabilistic potential theory, the Feynman--Kac formula, and degenerate elliptic equations.

\end{abstract}
\maketitle
\enlargethispage{2\baselineskip}
 \tableofcontents
\section{Introduction}

This paper studies the behavior of time fractional evolution equations as the order of the Caputo time derivative tends to zero: the solutions approach a resolvent value, and we determine both this resolvent limit and its stochastic interpretation, clarifying the connections with degenerate elliptic equations, probabilistic potential theory and the theory of stable subordinators. The scalar prototype of the problem is a classical asymptotic of Mittag-Leffler and Hardy, which we upgrade to an expansion of the Mittag-Leffler function to all orders in its order.

Let $X$ be a complex Banach space and let $A$ be a closed, densely defined linear operator on $X$. Consider the fractional Cauchy problem
\begin{equation}\label{FCP}
\begin{cases}
D_t^\alpha u(t)=Au(t),\quad t>0,\\
u(0)=\phi,
\end{cases}
\end{equation}
where $D_t^\alpha$, $0<\alpha<1$, is the Caputo fractional derivative (see Section \ref{sec-prelim} for the definition). Equations of this type are standard models for subdiffusion and for transport with memory, and they describe the scaling limits of continuous-time random walks; see e.g. \cite{Meerschaert2012book}. When $A$ generates a $C_0$-semigroup $(T(t))_{t\geq0}$, the solution of (\ref{FCP}) admits, under suitable assumptions, the subordination representation
$$
u_\alpha(t)=\int_0^\infty p_\alpha(s,t)T(s)\phi\,ds,
$$
where $p_\alpha(\cdot,t)$ is the probability density of an inverse $\alpha$-stable subordinator, see \cite{BM2001, BMN2009} and Section \ref{subsec-resolvent} below for the precise statement.

The question of how the solutions of (\ref{FCP}) depend on the order $\alpha$ has been studied intensively in the regime $\alpha\rightarrow 1^-$, where (\ref{FCP}) recovers the classical parabolic problem $u'=Au$: convergence of the Caputo derivative, of the Mittag-Leffler operators and of the solutions of linear and nonlinear fractional diffusion equations as $\alpha\rightarrow 1^-$ is by now well understood, see \cite{CT2023} and the references therein. At interior points of the parameter interval the solutions depend continuously on the order: if the corresponding problems are well-posed, then for every fixed initial datum and every fixed $t>0$, the solutions $u_\beta(t)$ converge to $u_{\alpha_0}(t)$ as $\beta\uparrow\alpha_0$ \cite{LiZheng2004}. The dependence on $\alpha$ is likewise central to inverse problems of recovering the order from additional data, see e.g. \cite{CNYY2009}. The opposite regime $\alpha\rightarrow 0^+$ has received much less attention, although related limiting densities were obtained in \cite{Orsingher2009}; see Remark \ref{rem-literature} below for this and related works. Its nature is essentially different. The point is that, unlike the Riemann--Liouville derivative, which satisfies ${}_{RL}D_t^\alpha\phi\rightarrow\phi$ as $\alpha\rightarrow 0^+$, the Caputo derivative does not converge to the identity operator. Indeed, since
\begin{equation*}
\widehat{D_t^\alpha\phi}(\lambda)=\lambda^\alpha\widehat{\phi}(\lambda)-\lambda^{\alpha-1}\phi(0)\longrightarrow \widehat{\phi}(\lambda)-\lambda^{-1}\phi(0) \qquad (\alpha\rightarrow 0^+),
\end{equation*}
one expects, and we shall prove,
\begin{equation}
\lim_{\alpha\rightarrow 0^+}D_t^\alpha\phi(t)=\phi(t)-\phi(0).
\end{equation}
The term $-\phi(0)$ records the collapse of the whole memory of the increments of $\phi$ into its total increment. As a consequence, the formal limit of (\ref{FCP}) as $\alpha\rightarrow 0^+$ is not $u=Au$ but the resolvent type equation
\begin{equation}\label{limit-eq-intro}
u(t)-\phi=Au(t),
\end{equation}
so that, under appropriate boundedness assumptions and with $1\in\rho(A)$, the expected limit is
\begin{equation}\label{limit-sol-intro}
\lim_{\alpha\rightarrow0^+}u_\alpha(t)=(1-A)^{-1}\phi\qquad \text{for every } t>0,
\end{equation}
a limit which is independent of $t>0$. Note that $u_\alpha(0)=\phi$ for every $\alpha>0$, so that a boundary layer occurs at $t=0$.

The link between (\ref{limit-sol-intro}) and the Mittag-Leffler function is formal but direct: the Laplace identity
\begin{equation*}
\lambda^{\alpha-1}(\lambda^\alpha-A)^{-1}=\sum_{n=0}^\infty\lambda^{-\alpha n-1}A^n
\end{equation*}
Formally, termwise inversion gives
$$
u_\alpha(t)=E_\alpha(At^\alpha)\phi
=\sum_{n=0}^\infty t^{\alpha n}A^n\phi/\Gamma(1+\alpha n).
$$
The resolvent family calculus of Sections \ref{subsec-resolvent} and \ref{subsec-rates} makes this computation rigorous. Taking the limit term by term yields $\sum_{n=0}^\infty A^n\phi=(1-A)^{-1}\phi$. The asymptotic behavior of the solutions as $\alpha\rightarrow0^+$ therefore reduces to that of $E_\alpha(z)$ in its order. In particular, the limit (\ref{limit-sol-intro}) rests on $E_\alpha(z)\rightarrow(1-z)^{-1}$. It is instructive to place this limit in its historical perspective. The Mittag-Leffler function grew out of Mittag-Leffler's work on the analytic continuation of a monogenic function beyond the disc of convergence. His summation method assigns to a series $\sum a_n$ the value $\lim_{\alpha\rightarrow0^+}\sum a_n/\Gamma(1+\alpha n)$, and for the geometric series this reads
\begin{equation*}
E_\alpha(z)=\sum_{n=0}^\infty\frac{z^n}{\Gamma(1+\alpha n)}\longrightarrow\frac{1}{1-z}\qquad(\alpha\rightarrow0^+),
\end{equation*}
uniformly on compact subsets of $\mathbb C\setminus[1,\infty)$, the Mittag-Leffler star of $(1-z)^{-1}$; see \cite[\S4.11, pp.~77--79 and \S8.10, Theorem 135, pp.~190--191]{Hardy1949}. The degeneration $\alpha\rightarrow0^+$ of $E_\alpha$ is thus a classical theme. For every constant $\zeta$, the function $t\mapsto E_\alpha(\zeta t^\alpha)$ solves $D_t^\alpha u=\zeta u$ with $u(0)=1$, and is a basic building block of linear time fractional equations. Parameter derivatives at the endpoint are already present in \cite[Eq.~(99)]{Apelblat2020}, and the same formal all-orders coefficient expansion appears in \cite[art.~14, Eq.~(41)]{VanMieghem2021}; see Remark \ref{rem-ML-scope}. Here we establish finite-order expansions with controlled remainders on compact subsets of the summability star and uniform estimates on sectors. These estimates are the basis of the operator-valued rates and expansions below.

The purpose of this paper is to make the above heuristic precise and to explore its consequences. After the preliminaries in Section \ref{sec-prelim}, the main results are presented in the following three directions, corresponding respectively to Sections \ref{sec-ML}, \ref{sec-abstract} and \ref{sec-concrete}.

(i) We establish a full expansion of the Mittag-Leffler function in its order at the endpoint $\alpha=0$, with rigorous finite-order remainder estimates on compact subsets of its summability star and uniform estimates on sectors (Theorem \ref{thm-ML-expansion}). Specializing to the scalar solution curve $E_\alpha(-\xi t^\alpha)$ yields the corresponding all-orders expansion and, in particular, a sharp first-order rate (Proposition \ref{prop-rate}). This scalar expansion provides the technical engine for the subsequent analysis.

(ii) For the abstract Caputo problem, we identify the limit equation as $u-\phi=Au$ and prove that $u_\alpha(t)\to(1-A)^{-1}\phi$ for every $t>0$, uniformly on compact time intervals (Theorem \ref{thm-homogeneous}). A uniform time bound on the solutions for one order suffices for this conclusion. The same boundedness assumption also forces sectoriality of $A$, which yields operator-norm convergence with rate $O(\alpha)$ and a full expansion in powers of $\alpha$, together with an inner expansion resolving the exponentially thin boundary layer at $t=0$ (Theorems \ref{thm-sectorial}, \ref{thm-unifbdd} and Proposition \ref{thm-inner}). The nonhomogeneous equation is treated as well (Theorem \ref{thm-nonhomogeneous}).

(iii) We provide a stochastic interpretation of the resolvent limit. The inverse $\alpha$-stable subordinator converges to a unit exponential random variable, and the limit solution is the $1$-potential of the underlying Markov process. This connects the limit equation with probabilistic potential theory, the Feynman--Kac formula and a half-space elliptic extension. Concrete realizations include the Laplacian, the fractional Poisson process and space-time fractional equations (Section \ref{sec-concrete}).

\section{Preliminaries}\label{sec-prelim}

Throughout the paper, $X$ is a complex Banach space (real spaces being complexified whenever contour methods are used), $A$ is a closed, densely defined linear operator on $X$ with domain $D(A)$, and $\rho(A)$ denotes the resolvent set of $A$. We write $\widehat{f}(\lambda)=\int_0^\infty e^{-\lambda t}f(t)\,dt$ for the (possibly vector-valued) Laplace transform.

\subsection{Fractional integrals and derivatives}

We begin with the definitions of fractional integrals and derivatives.

\begin{defn}
The Riemann--Liouville fractional integral of order $\alpha > 0$ is defined as
$$
J^\alpha_t f(t):= (g_\alpha* f)(t) = \int_0^t g_\alpha(t-s)f(s)\,ds,\quad\quad f\in L^1_{\mathrm{loc}}([0,+\infty);X),\quad t>0\ \text{a.e.}
$$
where
\begin{equation*}
g_\alpha(t)=\left\{
\begin{array}{ll}
\frac{ t^{\alpha-1}}{\Gamma(\alpha)},&t>0, \\
0,&t\leq 0.
\end{array}
\right.
\end{equation*}
Set moreover $J^0_t f(t)=f(t)$.
\end{defn}

\begin{defn}\label{def-Caputo}
The Caputo fractional derivative of order $\alpha > 0$ is defined as
$$
D^\alpha_t f(t):= J^{m-\alpha}_tD^m_tf(t), \quad t \in (0,T)\ \text{a.e.},\quad m-1<\alpha\leq m,
$$
if $f\in W^{m,1}((0,T);X)$, and the Riemann--Liouville fractional derivative of order $\alpha\in(0,1)$ is defined as
$$
{}_{RL}D^\alpha_t f(t):=\frac{d}{dt}J^{1-\alpha}_t f(t)
$$
with the derivative understood almost everywhere whenever $J^{1-\alpha}_t f$ is absolutely continuous. Moreover, we define
$$
{}_{RL}D^\alpha_t f(0):=\lim_{t \to 0} {}_{RL}D^\alpha_t f(t), \quad D^\alpha_t f(0):=\lim_{t \to 0} D^\alpha_t f(t)
$$
if the limits exist.
\end{defn}

If $0<\alpha<1$, the Laplace transforms of these derivatives read
\begin{equation*}
\widehat{D_t^\alpha u}(\lambda)=\lambda^\alpha \widehat{u}(\lambda)-\lambda^{\alpha-1}u(0),\qquad
\widehat{{}_{RL}D_t^\alpha u}(\lambda)=\lambda^\alpha \widehat{u}(\lambda)-J_t^{1-\alpha}u(0^+),
\end{equation*}
whenever the quantities involved make sense. In particular, ${}_{RL}D_t^\alpha$ and $D_t^\alpha$ differ by the way in which the initial value enters, and this difference is exactly what survives in the limit $\alpha\rightarrow 0^+$, see Section \ref{sec-Caputo-limit}.

\subsection{\texorpdfstring{$\alpha$}{alpha}-times resolvent families}\label{subsec-resolvent}

Consider the Volterra equation
\begin{equation}\label{Volterra}
u(t) = f(t) + \int_0^t g_\alpha(t-s) Au(s)\,ds, \quad t \ge 0
\end{equation}
where $f(t)$ is a continuous $X$-valued function.

\begin{defn} Let $u(t) : \mathbb R_+ \to X$ be continuous.

(1) $u(t)$ is called a strong solution of (\ref{Volterra}) if $u(t) \in D(A)$  and (\ref{Volterra}) holds for $t \ge 0$;

(2) $u(t)$ is called a mild solution of (\ref{Volterra}) if $(g_{\alpha} *u)(t) \in D(A)$ and $u(t) = f(t) + A(g_{\alpha} *u)(t) $ for $t \ge 0$.
\end{defn}

The notion of solution family for (\ref{Volterra}) is taken from \cite{Baj, Pruss}.
\begin{defn}
A family $\{S_\alpha(t)\}_{t \ge 0}\subset B(X)$ is called an $\alpha$-times resolvent family for the operator $A$ (or generated by $A$) if the following conditions are satisfied:

$(1)$ $S_\alpha(t)x: \mathbb R_+ \to X$ is continuous for every $x \in X$ and $S_\alpha(0)=I$;

$(2)$ $S_\alpha (t) D(A) \subset D(A)$ and $AS_\alpha(t) x = S_\alpha(t)Ax$ for all $x \in D(A)$ and $t \ge 0$;

$(3)$ the resolvent equation
\begin{equation*}
S_\alpha(t)x = x + (g_\alpha*S_\alpha)(t)Ax
\end{equation*}
holds for every $x \in D(A)$.
\end{defn}

\begin{rem}
(1) Since $A$ is closed and densely defined, it is easy to show that for all $x \in X$, $(g_\alpha*S_\alpha)(t)x \in D(A)$ and
\begin{equation*}
S_\alpha(t)x = x + A(g_\alpha*S_\alpha)(t)x.
\end{equation*}

(2) When $S_\alpha(t)$ is exponentially bounded, i.e., $\| S_\alpha(t)\|\leq Me^{\omega t}$, we have
\begin{equation}\label{LT-resolventfamily}
\int_0^\infty e^{-\lambda t}S_\alpha(t)x\, dt=\lambda^{\alpha-1}(\lambda^\alpha-A)^{-1}x,\qquad \mathrm{Re}\, \lambda>\max\{\omega,0\},\quad x\in X,
\end{equation}
with the principal branch of $\lambda^\alpha$.
The Laplace transform of vector-valued functions is one of the main technical tools in Section \ref{sec-abstract} below; for its theory we refer to \cite[Chapter 1]{ABHN2011} and \cite[Section 0.2]{Pruss}.
\end{rem}

If $A$ generates an $\alpha$-times resolvent family and
$f\in W^{1,1}((0,T);X)$, then the unique mild solution of
(\ref{Volterra}) is
$$
u(t)=S_\alpha(t)f(0)
+\int_0^t S_\alpha(t-s)f'(s)\,ds
=\frac{d}{dt}\int_0^t S_\alpha(t-s)f(s)\,ds;
$$
see \cite[Proposition 1.2, p.~33]{Pruss}. The differentiability in this formula is
guaranteed by the stated regularity of $f$.
In particular, the mild solution to
\begin{equation}\label{ugalphaA}
u(t)= \phi + \int_0^t g_\alpha(t-s) Au(s)\,ds
\end{equation}
is given by $u(t)= S_\alpha(t)\phi$. In addition, if $\phi \in D(A)$, then $u(t)$ is also a strong solution. Applying $D_t^\alpha$ to (\ref{ugalphaA})---legitimate under the regularity requirements of the Caputo formulation---for $0<\alpha\leq2$, leads to the fractional equation of order $\alpha$
\begin{equation}\label{alphaequation}
\begin{cases}
D_t^{\alpha}u(t) = Au(t), \quad t>0\\
u(0)= \phi\quad (u_t(0)= 0 \mbox{ if } 1<\alpha \le 2).
\end{cases}
\end{equation}
The mild formulation of (\ref{alphaequation}) is (\ref{ugalphaA}); its well-posedness for initial values $\phi\in X$ is characterized by the existence of an $\alpha$-times resolvent family, and the unique mild solution is $u(t)=S_\alpha(t)\phi$. For details we refer to \cite{Baj}.

The passage from semigroups to fractional resolvent families is effected by subordination. If $A$ generates an exponentially bounded $\beta$-times resolvent family $S_\beta$ for some $\beta>0$, then for every $0<\alpha<\beta$ it also generates an $\alpha$-times resolvent family, given by the subordination formula
\begin{equation}\label{subordination-general}
S_\alpha(t)x=\int_0^\infty p_{\alpha/\beta}(s,t)\,S_\beta(s)x\,ds,\qquad t>0,\ x\in X,
\end{equation}
where $p_\nu(\cdot,t)$ is the probability density of the inverse $\nu$-stable subordinator $H_\nu(t)$ introduced in Section \ref{subsec-subordinator} below; see \cite[Theorem 3.1, p.~40]{Baj} and, for $\beta=1$, \cite{BM2001}. Although the cited statement is for $0<\beta\leq2$, the same Laplace-transform calculation, using (\ref{inverse-stable-pdf-Laplace}) with exponent $\alpha/\beta\in(0,1)$, verifies the formula for every $\beta>0$. The case $\beta=1$, in which $S_1(t)=T(t)$ is a $C_0$-semigroup,
\begin{equation}\label{subordination2}
S_\alpha(t)\phi=\int_0^\infty p_\alpha(s,t)\,T(s)\phi\,ds,
\end{equation}
is the one used throughout this paper. In particular, the generator of a bounded $C_0$-semigroup also generates a bounded $\alpha$-times resolvent family for every $\alpha\in(0,1)$, with $\|S_\alpha(t)\|\leq\sup_{s\geq0}\|T(s)\|$.

The scalar case $X=\mathbb{C}$, $A=-\xi$ of (\ref{alphaequation}) has the explicit solution $u(t)=E_\alpha(-\xi t^\alpha)\phi$, where $E_\alpha$ is the Mittag-Leffler function: being the scalar prototype of the $\alpha$-times resolvent families, it is the engine of the limit theory of this paper and is studied in detail in Section \ref{sec-ML}.

\subsection{Stable and inverse stable subordinators}\label{subsec-subordinator}

\begin{defn}
A {\it subordinator} $D(s)$ is a one-dimensional L{\'e}vy process taking values in $[0,\infty)$ that is nondecreasing (a.s.) and starts at the origin (a.s.).

\end{defn}

Let $0<\alpha<1$ and let $Y_\alpha (t)$ be an $\alpha$-stable subordinator, i.e.,
a one-dimensional L{\'e}vy process starting at $0$, with nondecreasing (a.s.) paths and Laplace exponent
$\lambda^\alpha$, $\mathbb{E} e^{-\lambda Y_\alpha(t)}= e^{-t \lambda^{\alpha}}$, see \cite{Applebaum2009}.

The inverse $\alpha$-stable subordinator $H_\alpha(t)=\inf\{x>0:Y_\alpha(x)>t\}$ is a continuous, nondecreasing and nonnegative process
with $H_\alpha(0)=0$. For $t>0$, denote the probability density function of $H_\alpha(t)$ by $p_\alpha(x,t)$; then, with $E_\alpha$ the Mittag-Leffler function defined in (\ref{ML-def}) below, we have
\begin{equation}\label{ML-Laplace}
\int_0^\infty e^{-\xi x}p_\alpha(x,t)\,dx= E_\alpha(-\xi t^\alpha)
\end{equation}
and
\begin{equation}\label{inverse-stable-pdf-Laplace}
\int_0^\infty e^{-\lambda t}p_\alpha(x,t)\,dt= \lambda^{\alpha -1} e^{-x \lambda^\alpha},
\end{equation}
see e.g. \cite{Meerschaert2012book}. In view of Lemma \ref{ML-limit} below, (\ref{ML-Laplace}) identifies the limit of the distributions of $H_\alpha(t)$ as $\alpha\rightarrow0^+$ with the unit exponential distribution; this is made precise in Proposition \ref{H-convergence}.

\section{The Mittag-Leffler function and its expansion in the order}\label{sec-ML}

This section studies the scalar side of the problem. The first subsection establishes the basic limit of the scalar solution, $E_\alpha(-\xi t^\alpha)\rightarrow(1+\xi)^{-1}$. The second subsection upgrades this classical limit into an expansion of $E_\alpha(z)$ in powers of the order, with quantitative remainder estimates (Theorem \ref{thm-ML-expansion}), and specializes it to the solution curve $z=-\xi t^\alpha$, yielding an all-orders expansion as well as the sharp first-order rate (Proposition \ref{prop-rate}). The section concludes with a vector-valued convergence of the subordination measures (Lemma \ref{lem-measure-weak}), which provides a common analytic basis for the abstract and stochastic applications.

\subsection{Definition and the basic limit}

The Mittag-Leffler function is the scalar prototype of the $\alpha$-times resolvent families: for $X=\mathbb{C}$ and $A=-\xi$, the solution of (\ref{alphaequation}) is $u(t)=E_\alpha(-\xi t^\alpha)\phi$, where
\begin{equation}\label{ML-def}
E_\alpha(z) = \sum_{n=0}^\infty \frac{z^n}{\Gamma(\alpha n +1)},\qquad z\in\mathbb{C},\ \alpha>0
\end{equation}
is the Mittag-Leffler function. The following elementary limit is basic for everything that follows: it is the scalar ancestor of the resolvent limit $(1-A)^{-1}$ of Theorem \ref{thm-homogeneous}. Note that $1/(1+\xi)$ is the Laplace transform of the unit exponential density $e^{-x}$.

\begin{lem}\label{ML-limit}
Let $t>0$ be fixed. Then for every real $\xi>-1$,
\begin{equation}\label{ML-limit-eq}
\lim_{\alpha\rightarrow 0^+}E_\alpha(-\xi t^\alpha)=\frac{1}{1+\xi}.
\end{equation}
More precisely, the convergence holds locally uniformly for real $\xi$ in $(-1,\infty)$, and indeed locally uniformly on the open right half-plane.
\end{lem}
\begin{proof}
The full statement follows directly from Hardy's theorem (Remark \ref{rem-ML-summability}): since $-\xi t^\alpha\rightarrow-\xi$ stays eventually in a compact subset of $\mathbb C\setminus[1,\infty)$, the uniform convergence there yields (\ref{ML-limit-eq}), locally uniformly in $\xi$ and uniformly for $t$ in compact subsets of $(0,\infty)$. However, to keep this basic limit self-contained, and because the domination argument below will be reused in Example \ref{ex-exp}, we give the following elementary proof.

Consider the $-1<\xi<1$ case.
 Since $t^\alpha\rightarrow 1$ as $\alpha\rightarrow0^+$, we can choose $\alpha_0>0$ such that $|\xi|\,t^{\alpha_0}<1$, and then $\rho:=|\xi|\max(1,t^{\alpha_0})<1$ because $|\xi|<1$. For all $\alpha\in(0,\alpha_0]$ and all $n\geq0$ we have $|\xi|^n t^{\alpha n}=(|\xi|\,t^\alpha)^n\leq \rho^n$: indeed, $t^\alpha\leq1$ if $t\leq1$, while $t^\alpha\leq t^{\alpha_0}$ if $t\geq1$. Moreover, $\Gamma(x)\geq\int_1^2e^{-t}\,dt=e^{-1}-e^{-2}$ for all $x\geq1$ (since $t^{x-1}\geq1$ for $t\geq1$), hence $\Gamma(\alpha n+1)\geq e^{-1}-e^{-2}$ and
$$
\left|\frac{(-\xi)^n t^{\alpha n}}{\Gamma(\alpha n+1)}\right|\leq \frac{\rho^n}{e^{-1}-e^{-2}},\qquad 0<\alpha\leq\alpha_0,\ n\geq0.
$$
By the Weierstrass $M$-test, the series (\ref{ML-def}) with $z=-\xi t^\alpha$ converges uniformly in $\alpha\in(0,\alpha_0]$. Therefore we can pass to the limit term by term: since $t^{\alpha n}\rightarrow1$ and $\Gamma(\alpha n+1)\rightarrow\Gamma(1)=1$ for each fixed $n$,
$$
\lim_{\alpha\rightarrow0^+}E_\alpha(-\xi t^\alpha)=\sum_{n=0}^\infty(-\xi)^n=\frac{1}{1+\xi}.
$$

Next consider the $\xi\geq1$ case.  By (\ref{ML-Laplace}), $E_\alpha(-\xi t^\alpha)$ is the Laplace transform of a probability density, hence analytic and bounded by $1$ in $\mathrm{Re}\,\xi>0$. By Step 1 it converges to $1/(1+\xi)$ on the interval $(0,1)$. Vitali's convergence theorem therefore yields locally uniform convergence on $\mathrm{Re}\,\xi>0$ to an analytic limit, which coincides with $1/(1+\xi)$ by the identity theorem. In particular (\ref{ML-limit-eq}) holds for every real $\xi>0$; together with Step 1, which covers $-1<\xi\leq0$, this proves the full statement.
\end{proof}

Written as $D_t^\alpha u=au$, the scalar equation is thus covered for every $a<1$ (with $a=-\xi$), with limit $(1-a)^{-1}$, whereas for $a\geq1$ the solution diverges (Remark \ref{rem-ML-summability} and the boundedness discussion following Remark \ref{rem-boundary-layer}).

\begin{rem}\label{rem-ML-summability}
Applied to the geometric series, Mittag-Leffler's summation method
gives $E_\alpha(z)\to(1-z)^{-1}$ uniformly on compact subsets of
$\mathbb C\setminus[1,\infty)$; see
\cite[\S4.11, p.~79; \S8.10, Theorem~135]{Hardy1949} and the
historical account in \cite{GKMR2020}. The cut cannot be included:
for real $x\geq1$ and every $N$,
$$
E_\alpha(x)\geq\sum_{n=0}^N\frac{x^n}{\Gamma(1+\alpha n)}
\longrightarrow\sum_{n=0}^N x^n,
$$
so $E_\alpha(x)\to+\infty$.
\end{rem}

\subsection{Expansion in the order}\label{subsec-ML-expansion}

Lemma \ref{ML-limit} identifies the basic limit $E_\alpha(z)\rightarrow(1-z)^{-1}$ as $\alpha\rightarrow0^+$. For the operator-valued applications of Section \ref{subsec-rates} we need quantitative versions of this limit, in the form of an expansion of $E_\alpha(z)$ in powers of the order $\alpha$ with a controlled remainder, valid outside a neighborhood of $[1,\infty)$. We state it for general $z\in\mathbb{C}\setminus[1,\infty)$, then the scalar relaxation case corresponds to $z=-\xi t^\alpha$, $\xi>0$.

The coefficients are built from the Taylor coefficients $c_j$ of $1/\Gamma(1+x)$ at $x=0$. The reciprocal gamma function $1/\Gamma(s)$ is an entire function of $s\in\mathbb C$, with simple zeros at $s=0,-1,-2,\dots$ and no others \cite[\S5.2(i)]{NIST}. In particular $1/\Gamma(1+x)$ is entire in $x$, so its Maclaurin series converges and represents it for every $x\in\mathbb C$. Explicitly, the digamma function $\psi=\Gamma'/\Gamma$ admits the Maclaurin series
$$
\psi(1+x)=-\gamma+\sum_{j=2}^\infty(-1)^j\,\zeta(j)\,x^{j-1},\qquad |x|<1,
$$
where $\gamma=-\Gamma'(1)$ is the Euler constant and $\zeta$ denotes the Riemann zeta function, see \cite[(5.7.4)]{NIST}. Since $\psi(1+x)=\frac{d}{dx}\log\Gamma(1+x)$ and $\log\Gamma(1)=0$, termwise integration followed by exponentiation gives
\begin{equation}\label{Gamma-reciprocal}
\frac{1}{\Gamma(1+x)}=\sum_{j=0}^\infty c_jx^j=\exp\Big(\gamma x-\sum_{j=2}^\infty(-1)^j\frac{\zeta(j)}{j}\,x^j\Big),\qquad |x|<1,
\end{equation}
using $\zeta(2)=\pi^2/6$, $c_0=1$, $c_1=\gamma$, $c_2=\frac{\gamma^2}{2}-\frac{\pi^2}{12}$, $c_3=\frac{\gamma^3}{6}-\frac{\gamma\pi^2}{12}+\frac{\zeta(3)}{3}$, \dots. A self-contained derivation of the same series, including a closed form for the coefficients $c_j$ in terms of the Stirling numbers of the first kind, can be found in \cite{Ulgenes2025}. The restriction $|x|<1$ in (\ref{Gamma-reciprocal}) comes from the exponential form on its right, $\log\Gamma(1+x)$ having its nearest singularity at $x=-1$; the power series itself converges for every $x\in\mathbb C$, as $1/\Gamma(1+x)$ is entire. The appearance of the reciprocal gamma function at this point is the exact counterpart of the classical theory: in Hardy's proof of Mittag-Leffler summability ``the details naturally depend on the asymptotic properties of the gamma-function of a complex variable'' \cite[\S4.11, p.~79]{Hardy1949}; Theorem \ref{thm-ML-expansion} makes this dependence on the gamma function explicit and quantitative through the series (\ref{Gamma-reciprocal}). The other ingredient is the polylogarithms of negative order
$$
\mathrm{Li}_{-k}(z):=\sum_{n=1}^\infty n^kz^n=\dfrac{zQ_k(z)}{(1-z)^{k+1}},
$$
where, for $k\geq1$, $Q_k$ is the $k$-th Eulerian polynomial. Indeed, $\sum_{n\geq1}n^kz^n=(z\frac{d}{dz})^k\frac{1}{1-z}$, and an induction on $k$ gives $Q_1=1$ and $Q_{k+1}(z)=(1+kz)Q_k(z)+z(1-z)Q_k'(z)$. Thus $zQ_k(z)$ is the Eulerian polynomial in the normalization of \cite[\S1.4, Proposition~1.4.4]{StanleyEC1}; its coefficients count permutations by descents and are in particular nonnegative {\rm(}$Q_1=1$, $Q_2(z)=1+z$, $Q_3(z)=1+4z+z^2$, \dots{\rm)}. The series converges for $|z|<1$, and the rational expression on the right, meaningful for all $z\neq1$, gives the analytic continuation by which $\mathrm{Li}_{-k}(z)$ is understood for $z\in\mathbb{C}\setminus[1,\infty)$ below. For later use, set, with the coefficients $c_j$ of (\ref{Gamma-reciprocal}),
$$
P_k(u):=\sum_{j=0}^k\frac{c_j}{(k-j)!}\,u^{k-j},\qquad P_0=1,\ \ P_1(u)=u+\gamma,\ \ P_2(u)=\frac{u^2}{2}+\gamma u+\frac{\gamma^2}{2}-\frac{\pi^2}{12},\ \dots
$$
Heuristically, we then have
\begin{align*}
E_\alpha(z)&=\sum_{n=0}^\infty \frac{z^n}{\Gamma(1+\alpha n)}=\sum_{n=0}^\infty z^n \sum_{k=0}^\infty c_k (\alpha n)^k
= \sum_{k=0}^\infty c_k \alpha^k \sum_{n=0}^\infty z^n n^k\\
&=c_0 \sum_{n=0}^\infty z^n +\sum_{k=1}^\infty c_k \alpha^k \sum_{n=0}^\infty z^n n^k=\frac{c_0}{1-z}+\sum_{k=1}^\infty c_k \alpha^k \frac{zQ_k(z)}{(1-z)^{k+1}}\\
&=\frac{1}{1-z}+\sum_{k=1}^K\alpha^k\,c_k\,\mathrm{Li}_{-k}(z)+\sum_{k=K+1}^\infty c_k \alpha^k \frac{zQ_k(z)}{(1-z)^{k+1}}.
\end{align*}
Here the substitution of the Maclaurin series of $1/\Gamma(1+x)$ at $x=\alpha n$ is legitimate for every $n$, since $1/\Gamma(1+x)$ is entire, as noted above. The interchange of the two summations, however, is purely formal: the resulting series is asymptotic at best. Theorem \ref{thm-ML-expansion} below turns the truncation into a rigorous statement.

\begin{lem}\label{lem-ML-Hankel}
Let $\theta=3\pi/4$ and let $C$ consist of the arc $\{|\zeta|=1,\ |\arg\zeta|\leq\theta\}$ and the two rays $\{\arg\zeta=\pm\theta,\ |\zeta|\geq1\}$, oriented so that the region containing $0$ and the negative real axis lies on the left. Then Hankel's formula
\begin{equation}\label{R-Gamma-Hankel}
\frac{1}{2\pi i}\int_C e^\zeta\,\zeta^{-s}\,d\zeta=\frac{1}{\Gamma(s)},\qquad s>0,
\end{equation}
holds, and for every $0<\alpha<1$ the Mittag-Leffler function admits the contour representation
\begin{equation}\label{ML-Hankel}
E_\alpha(z)=\frac{1}{2\pi i}\int_C \frac{e^\zeta\,\zeta^{\alpha-1}}{\zeta^\alpha-z}\,d\zeta,
\end{equation}
valid for $\{|z|<1\}\ \cup\ \{|\arg z|\geq\delta\}$ where $\alpha\theta<\delta<\pi$.
\end{lem}
\begin{proof}
The first identity is Hankel's formula for the reciprocal gamma function
\cite[(5.9.2)]{NIST}, with its rays deformed to arguments $\pm\theta$.
The deformation is valid because $e^\zeta$ decays exponentially on the
intervening rays. The integral converges locally uniformly for
$s\in\mathbb C$, so the formula also holds for every complex $s$.

For $|z|<1$, expand
$(\zeta^\alpha-z)^{-1}=\sum_{n\geq0}z^n\zeta^{-\alpha(n+1)}$ on $C$.
The expansion is uniformly absolutely convergent there, and its product
with $e^\zeta\zeta^{\alpha-1}$ has an integrable majorant. Termwise
integration and Hankel's formula give (\ref{ML-def}).

It remains to extend (\ref{ML-Hankel}) beyond the disk. Fix
$0<\alpha<1$ and set
$$
D_\alpha:=\{|z|<1\}\ \cup\ \{|\arg z|>\alpha\theta\}.
$$
The disk and the open sector overlap (both contain $-1/2$, say), so
$D_\alpha$ is a connected open set. We claim that the right-hand side
of (\ref{ML-Hankel}) defines a holomorphic function on $D_\alpha$.
First, the denominator never vanishes there: since $|\zeta|\geq1$ on
$C$ and $|\arg\zeta^\alpha|=\alpha|\arg\zeta|\leq\alpha\theta$
(principal value; no branch ambiguity arises since
$\theta=3\pi/4<\pi$), every value taken by $\zeta^\alpha$ on $C$
satisfies $|\zeta^\alpha|\geq1$ and $|\arg\zeta^\alpha|\leq\alpha\theta$,
whereas a point $z\in D_\alpha$ satisfies either $|z|<1$ or
$|\arg z|>\alpha\theta$; hence $\zeta^\alpha-z\neq0$. Second, the
integral converges locally uniformly in $z$: if $K\subset D_\alpha$
is compact, then $|\zeta^\alpha-z|\geq c_K>0$ on the bounded part of
$C$ by continuity, and
$|\zeta^\alpha-z|\geq\tfrac12|\zeta|^\alpha$ for $|\zeta|$
sufficiently large, uniformly for $z\in K$; hence the integrand is
dominated by the $z$-independent integrable majorant
$C_K\,e^{\mathrm{Re}\,\zeta}(1+|\zeta|)^{-1}$ on $C$, and holomorphy
follows by differentiation under the integral sign. Since $E_\alpha$
is entire and agrees with this function on the disk $\{|z|<1\}$ by the
series argument above, the identity theorem extends (\ref{ML-Hankel})
to all of $D_\alpha$. Finally, whenever $\alpha\theta<\delta<\pi$ one
has
$\{|\arg z|\geq\delta\}\subset\{|\arg z|>\alpha\theta\}\subset D_\alpha$,
which completes the proof.
\end{proof}

\begin{thm}\label{thm-ML-expansion}
For every integer $K\geq0$ (with the convention that a sum over the empty index range, arising for $K=0$, equals zero), the Mittag-Leffler function admits the following expansion as $\alpha\rightarrow0^+$:
\begin{equation}\label{ML-all-orders-eq}
E_\alpha(z)=\frac{1}{1-z}+\sum_{k=1}^K\alpha^k\,c_k\,\mathrm{Li}_{-k}(z)+R_K(z,\alpha),\qquad \alpha>0,\ z\in\mathbb{C}\setminus[1,\infty),
\end{equation}
with remainder $R_K(z,\alpha)=O(\alpha^{K+1})$ in the following precise sense.

{\rm (i)} For every $\rho\in(0,1)$ there exists $C_{K,\rho}>0$ such that
$$
|R_K(z,\alpha)|\leq C_{K,\rho}\,\alpha^{K+1},\qquad |z|\leq\rho,\ 0<\alpha\leq1.
$$

{\rm (ii)} For every $\delta\in(0,\pi)$ there exist $\alpha_0(\delta)>0$ and $C_{K,\delta}>0$ such that, for all $0<\alpha\leq\alpha_0(\delta)$ and all $z$ with $|\arg z|\geq\delta$,
$$
|E_\alpha(z)|\leq\frac{C_\delta}{1+|z|},\qquad
|R_K(z,\alpha)|\leq\frac{C_{K,\delta}\,\alpha^{K+1}\,|z|\,(1+|z|)^{K-1}}{|1-z|^{K+1}}.
$$

{\rm (iii)} The expansion {\rm(\ref{ML-all-orders-eq})} holds uniformly, with rate $O(\alpha^{K+1})$, on every compact subset of $\mathbb{C}\setminus[1,\infty)$.

{\rm (iv)}  The sector estimate of {\rm(ii)} is stable under differentiation in $z$: for $j=1,2$, for every $\delta\in(0,\pi)$ there exist $\alpha_0(\delta)>0$ and $C_{K,\delta,j}>0$ such that
$$
\bigl|\partial_z^jR_K(z,\alpha)\bigr|\leq\frac{C_{K,\delta,j}\,\alpha^{K+1}\,(1+|z|)^{K-1+j}}{|1-z|^{K+1+j}},\qquad |\arg z|\geq\delta,\ 0<\alpha\leq\alpha_0(\delta).
$$
\end{thm}

\begin{proof}
(i) We first claim that
\begin{equation}\label{rK(x)-difference}
\left|\frac1{\Gamma(1+x)}-\sum_{k=0}^K c_kx^k\right|
\leq C_Kx^{K+1},\qquad x\geq0.
\end{equation}
For $0\leq x\leq1$ this is Taylor's theorem with the Lagrange
remainder, $1/\Gamma(1+x)$ being smooth. For $x\geq1$ the two terms
are estimated separately. First, $1/\Gamma(1+x)$ is bounded on
$[0,\infty)$, being continuous with $\Gamma(1+x)\rightarrow+\infty$ as
$x\rightarrow+\infty$ (a crude explicit bound is
$\Gamma(1+x)\geq e^{-1}-e^{-2}$, cf.\ the proof of Lemma
\ref{ML-limit}); hence $|1/\Gamma(1+x)|\leq C_Kx^{K+1}$ for $x\geq1$.
Second, $|\sum_{k=0}^Kc_kx^k|\leq C_Kx^K\leq C_Kx^{K+1}$ for $x\geq1$. The
claim follows for all $x\geq0$.

Write $r_K(x)$ for the difference in the claim \eqref{rK(x)-difference} and substitute
$x=\alpha n$ into (\ref{ML-def}):
$$
E_\alpha(z)=\sum_{n=0}^\infty z^n
\Bigl(\sum_{k=0}^Kc_k(\alpha n)^k+r_K(\alpha n)\Bigr).
$$
For $|z|\leq\rho$ each of the $K+1$ subseries converges absolutely,
since $\sum_{n=0}^\infty n^k\rho^n<\infty$ for each $k\leq K$. The
sum over the finite index set $\{0,\dots,K\}$ therefore interchanges
with the sum over $n$. The
main terms are exactly those of (\ref{ML-all-orders-eq}): for
$k\geq1$, $\sum_{n=0}^\infty n^kz^n=\mathrm{Li}_{-k}(z)$, and the
$k=0$ term is $c_0\sum_{n\geq0}z^n=c_0/(1-z)=1/(1-z)$ since $c_0=1$.
The remainder defined by subtraction in (\ref{ML-all-orders-eq}) is
thus $R_K(z,\alpha)=\sum_{n=0}^\infty z^n r_K(\alpha n)$, the series
converging absolutely by the estimate that follows, and
$$
|R_K(z,\alpha)|\leq C_K\alpha^{K+1}
\sum_{n=0}^\infty n^{K+1}|z|^n
\leq C_K\alpha^{K+1}\sum_{n=0}^\infty n^{K+1}\rho^n
=:C_{K,\rho}\,\alpha^{K+1},
$$
which proves (i).

(ii) Use the contour $C$ of Lemma \ref{lem-ML-Hankel}, put
$\alpha_0(\delta)=\delta/(2\theta)<1$, and define, for
$0\leq\beta\leq\alpha_0(\delta)$ (including $\beta=0$) and
$|\arg z|\geq\delta$,
\begin{equation*}
F(\beta,z):=\frac{1}{2\pi i}\int_C e^\zeta\,\zeta^{-1}\,
\frac{\zeta^\beta}{\zeta^\beta-z}\,d\zeta.
\end{equation*}
By Lemma \ref{lem-ML-Hankel}, $F(\alpha,z)=E_\alpha(z)$ for
$0<\alpha\leq\alpha_0(\delta)$, since $\alpha_0(\delta)<\delta/\theta$;
Hankel's formula \eqref{R-Gamma-Hankel}  at $s=1$ gives $F(0,z)=(1-z)^{-1}$.

For $\zeta\in C$ and $0\leq\beta\leq\alpha_0(\delta)$ one has
$|\arg\zeta^\beta|=\beta|\arg\zeta|\leq\beta\theta\leq\delta/2$
(principal value; no branch ambiguity arises since
$\theta=3\pi/4<\pi$), so the angle between $\zeta^\beta$ and $z$ is at
least $\delta/2$ whenever $|\arg z|\geq\delta$ (see Figure
\ref{fig-contour}). Two nonzero numbers
separated by an angle $\varphi$ satisfy
$|u-v|\geq\sin(\varphi/2)(|u|+|v|)$, as follows from
$|u-v|^2=(|u|+|v|)^2-4|u||v|\cos^2(\varphi/2)$. Hence
$$
|\zeta^\beta-z|\geq\sin(\delta/4)
\bigl(|\zeta|^\beta+|z|\bigr).
$$
Equivalently, $w:=z\zeta^{-\beta}$ makes an angle at least $\delta/2$
with $1$, and
\begin{equation}\label{angle-gap-w}
|1-w|\geq\sin(\delta/4)\bigl(1+|w|\bigr).
\end{equation}
Since $|\zeta|\geq1$ on $C$, we have $|\zeta|^\beta\leq|\zeta|^{\alpha_0}$. If $|z|\leq1$, then $1+|z|\leq2$, so the quotient $|\zeta|^\beta/(|\zeta|^\beta+|z|)$ is at most $1\leq2|\zeta|^{\alpha_0}/(1+|z|)$. If $|z|\geq1$, then $1+|z|\leq2|z|$, so the quotient is at most $|\zeta|^\beta/|z|\leq2|\zeta|^{\alpha_0}/(1+|z|)$. Thus
$$
\left|\frac{\zeta^\beta}{\zeta^\beta-z}\right|
\leq\frac{C_\delta|\zeta|^{\alpha_0}}{1+|z|}.
$$
Multiplying by $e^\zeta\zeta^{-1}$ and integrating over $C$ proves the
bound for $E_\alpha(z)$, since
$$
\int_C e^{\mathrm{Re}\,\zeta}|\zeta|^{\alpha_0-1}\,|d\zeta|<\infty.
$$

\begin{figure}[t]
\centering
\includegraphics[width=0.66\textwidth]{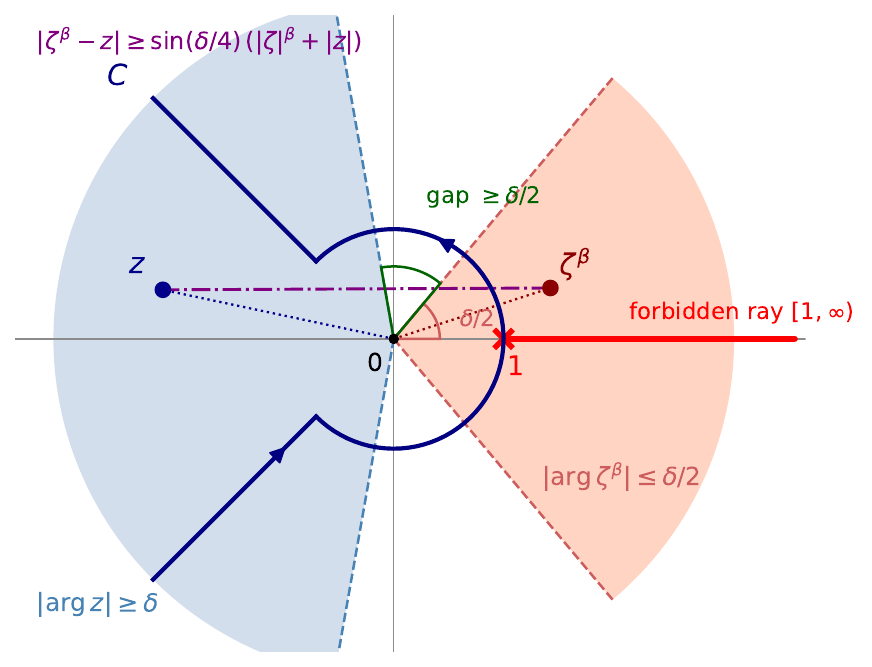}
\caption{The geometry of the proof of Theorem \ref{thm-ML-expansion}(ii). The contour $C$ consists of the arc $\{|\zeta|=1,\ |\arg\zeta|\leq3\pi/4\}$ and the rays $\arg\zeta=\pm3\pi/4$, traversed from the lower ray to the upper one. For $\zeta\in C$ and $0\leq\beta\leq\alpha_0(\delta)$, the point $\zeta^\beta$ stays in the narrow sector $|\arg\cdot|\leq\delta/2$ around the positive real axis, while $z$ lies in the wide sector $|\arg z|\geq\delta$ around the negative real axis; the angle between $\zeta^\beta$ and $z$ is therefore at least $\delta/2$, which yields the lower bound $|\zeta^\beta-z|\geq\sin(\delta/4)\,(|\zeta|^\beta+|z|)$. The bold red ray $[1,\infty)$ is the forbidden cut of the Mittag-Leffler limit.}
\label{fig-contour}
\end{figure}

We turn to the expansion. Writing $w=z\zeta^{-\beta}$, one has
$dw/d\beta=(-\log\zeta)\,w$, hence
$d/d\beta=(-\log\zeta)(w\,d/dw)$ and, for every $m\geq1$,
\[
\frac{d^m}{d\beta^m}
\frac{\zeta^\beta}{\zeta^\beta-z}=\frac{d^m}{d\beta^m}\frac{1}{1-z\zeta^{-\beta}}
=(-\log\zeta)^m\Bigl(w\frac{d}{dw}\Bigr)^m\frac{1}{1-w}
=(-\log\zeta)^m\mathrm{Li}_{-m}(z\zeta^{-\beta}),
\]
where the last equality, valid a priori for $|w|<1$ as a series
identity, persists for all $w\neq1$ since both sides are rational in
$w$. Since $Q_m$ has degree $m-1$, (\ref{angle-gap-w}) gives
$$
|\mathrm{Li}_{-m}(w)|
=\frac{|w|\,|Q_m(w)|}{|1-w|^{m+1}}\leq \frac{C'|w|(1+|w|)^{m-1}}{|1-w|^{m+1}}
\leq \frac{C_{m,\delta}|w|}{(1+|w|)^2}
=\frac{C_{m,\delta}|z|\,|\zeta|^\beta}{(|\zeta|^\beta+|z|)^2}
\leq \frac{C_{m,\delta}|\zeta|^{\alpha_0}|z|}{(1+|z|)^2},
$$
uniformly for $0\leq\beta\leq\alpha_0(\delta)$ and $|\arg z|\geq\delta$.

We now Taylor-expand $F$ in $\beta$ at $\beta=0$. Integrating the derivative identity above against the dominating factors shows that every derivative
$$
\frac{\partial^m F}{\partial\beta^m}(\beta,z)
=\frac{1}{2\pi i}\int_C e^\zeta\zeta^{-1}(-\log\zeta)^m
\mathrm{Li}_{-m}(z\zeta^{-\beta})\,d\zeta,\qquad 1\leq m\leq K+1,
$$
satisfies the uniform bound
\begin{equation}\label{F-deriv-bound}
\Bigl|\frac{\partial^m F}{\partial\beta^m}(\beta,z)\Bigr|
\leq C_{m,\delta}\,\frac{|z|}{(1+|z|)^2},
\qquad 0\leq\beta\leq\alpha_0(\delta),\ \ |\arg z|\geq\delta,
\end{equation}
the extra $\zeta$-dependent factors being absorbed by
$$
\int_C e^{\mathrm{Re}\,\zeta}|\zeta|^{\alpha_0-1}
\bigl(1+|\log\zeta|^{K+1}\bigr)\,|d\zeta|<\infty.
$$
Hence $F$ is $C^{K+1}$ in $\beta$ on $[0,\alpha_0(\delta)]$, and
Taylor's formula with integral remainder applies:
\begin{equation}\label{F-Taylor}
F(\alpha,z)=\sum_{k=0}^K\frac{\alpha^k}{k!}\,
\partial_\beta^kF(0,z)
+\frac{1}{K!}\int_0^\alpha(\alpha-\beta)^K\,
\partial_\beta^{K+1}F(\beta,z)\,d\beta.
\end{equation}

It remains to evaluate the coefficients at $\beta=0$. There the
factor $\mathrm{Li}_{-k}(z\zeta^{-\beta})=\mathrm{Li}_{-k}(z)$ is
independent of $\zeta$ and comes out of the integral. To differentiate
Hankel's formula \eqref{R-Gamma-Hankel} under the integral sign, observe
that $|\zeta^{-s}|\leq|\zeta|^{\alpha_0/2-1}$ on $C$ whenever
$\operatorname{Re}s\geq1-\alpha_0/2$. Since
$e^{\operatorname{Re}\zeta}|\zeta|^{\alpha_0/2-1}
(1+|\log\zeta|)^k$ is integrable on $C$, differentiation is legitimate
in a neighbourhood of $s=1$ and gives
$$
\frac1{2\pi i}\int_C e^\zeta\zeta^{-1}
(-\log\zeta)^k\,d\zeta
=\left.\frac{d^k}{ds^k}\frac1{\Gamma(s)}\right|_{s=1}
=k!c_k,
$$
so that $\partial_\beta^kF(0,z)=k!\,c_k\,\mathrm{Li}_{-k}(z)$ for
$k\geq1$, while $F(0,z)=(1-z)^{-1}$ as already noted. Hence
(\ref{F-Taylor}) is exactly the identity (\ref{ML-all-orders-eq}), the
remainder defined there being
$$
R_K(z,\alpha)=\frac{1}{K!}\int_0^\alpha(\alpha-\beta)^K
\frac{\partial^{K+1}F}{\partial\beta^{K+1}}(\beta,z)\,d\beta.
$$
The case $m=K+1$ of (\ref{F-deriv-bound}) then gives the equivalent
sector bound
\begin{equation}\label{ML-sector-remainder}
|R_K(z,\alpha)|\leq C_{K,\delta}\alpha^{K+1}
\frac{|z|}{(1+|z|)^2},
\end{equation}
the factor $\alpha^{K+1}/(K+1)!$ from the $\beta$-integral being
absorbed into $C_{K,\delta}$. Since $|1-z|^{K+1}\leq(1+|z|)^{K+1}$, one has
$|z|/(1+|z|)^2\leq|z|(1+|z|)^{K-1}/|1-z|^{K+1}$, and the stated
estimate in (ii) follows.

(iii) Every compact subset of $\mathbb C\setminus[1,\infty)$ is
contained in the union of a disk $\{|z|\leq\rho\}$, with $\rho<1$,
and a sector $\{|\arg z|\geq\delta\}$, with $\delta>0$. Explicitly,
let $d>0$ be its distance to $[1,\infty)$, let $L\geq1$ be an upper
bound for the moduli of its points, and put $q=\min(1,d)$,
$\rho=1-q/2$ and $\delta=\arcsin(q/(4L))$. Indeed, if
$|z|>\rho$ and $|\arg z|<\delta$, set $r=|z|\leq L$. Then
$\max(1,r)\in[1,\infty)$ and
$$
\operatorname{dist}(z,[1,\infty))
\leq |z-r|+(1-r)_+
<2L\sin(\delta/2)+\frac q2
\leq L\delta+\frac q2
\leq\left(\frac\pi8+\frac12\right)q<d,
$$
where we used $\delta\leq(\pi/2)\sin\delta=\pi q/(8L)$. This
contradiction proves the covering assertion. Parts (i) and (ii)
therefore give the claim, the $z$-dependent factor in the bound of
(ii) being bounded on the compact set.

(iv) The remainder $R_K(\cdot,\alpha)$ is holomorphic on
$\mathbb C\setminus[1,\infty)$, since $E_\alpha$ is entire and each term
of (\ref{ML-all-orders-eq}) is rational with its only pole at
$z=1$. For $|z|\leq1/2$, Cauchy's estimates on the circle of radius
$1/4$ about $z$, contained in $\{|w|\leq3/4\}$, combined with (i) for
$\rho=3/4$, give
$|\partial_z^jR_K(z,\alpha)|\leq j!\,4^jC_{K,3/4}\,\alpha^{K+1}$; this
implies the claimed bound on $|z|\leq1/2$, where $1+|z|$ and $|1-z|$
are bounded from above and from below by positive constants.

For $|z|>1/2$ with $|\arg z|\geq\delta$, choose
$0<\eta_\delta<\tfrac12\sin(\delta/2)$: the circle of radius
$\eta_\delta|z|$ about $z$ then lies in $\{|\arg w|\geq\delta/2\}$,
and on it $|z|/2\leq|w|\leq3|z|/2$. Cauchy's estimates and
(\ref{ML-sector-remainder}), applied with the angle $\delta/2$ (which
requires $\alpha\leq\alpha_0(\delta/2)$ and accounts for the decrease
of $\alpha_0(\delta)$), give
$$
|\partial_z^jR_K(z,\alpha)|
\leq\frac{j!}{(\eta_\delta|z|)^j}\,
C_{K,\delta/2}\,\alpha^{K+1}\sup_{|w-z|=\eta_\delta|z|}
\frac{|w|}{(1+|w|)^2}
\leq C_{K,\delta,j}'\,\alpha^{K+1}\,
\frac{|z|^{1-j}}{(1+|z|)^2},
$$
since $|w|/(1+|w|)^2\leq C|z|/(1+|z|)^2$ when $|w|$ is comparable to
$|z|$. For $j=1$ the last factor is $(1+|z|)^{-2}$; for $j=2$ one uses
$|z|^{-1}(1+|z|)\leq3$, valid for $|z|>1/2$. In both cases
$$
|\partial_z^jR_K(z,\alpha)|\leq
\frac{C_{K,\delta,j}\alpha^{K+1}}{(1+|z|)^{j+1}},
\qquad j=1,2.
$$
Finally,
$|1-z|^{K+1+j}\leq(1+|z|)^{K+1+j}\leq(1+|z|)^{K+2j}$ gives
$1/(1+|z|)^{j+1}\leq(1+|z|)^{K-1+j}/|1-z|^{K+1+j}$, and (iv) follows.
\end{proof}

\begin{rem}\label{rem-sector-equivalent}
For $|\arg z|\geq\delta$, the angle-gap inequality gives
$\sin(\delta/2)(1+|z|)\leq|1-z|\leq1+|z|$.
Thus the remainder bound in Theorem \ref{thm-ML-expansion}(ii) is
equivalent to
$$
|R_K(z,\alpha)|\leq
C_{K,\delta}\alpha^{K+1}\frac{|z|}{(1+|z|)^2},
$$
the form used in the operator arguments below. Its decay at infinity
is sharp: applying the theorem with $K=2$ and using
$\mathrm{Li}_{-2}(-r)=r(r-1)/(1+r)^3$ gives
$$
R_1(-r,\alpha)=\frac{c_2\alpha^2}{r}
+O\left(\frac{\alpha^2}{r^2}+\frac{\alpha^3}{r}\right),
\qquad r\geq1,
$$
where $c_2\neq0$. The bounds in part (iv) are equivalently
$|\partial_z^jR_K(z,\alpha)|\leq
C_{K,\delta,j}\alpha^{K+1}/(1+|z|)^2$ for $j=1,2$.

The domains of parts (i) and (ii) of Theorem \ref{thm-ML-expansion}
overlap: a point with $|z|\leq\rho<1$ and $|\arg z|\geq\delta$ is
covered by both. The two estimates serve complementary purposes: (i)
gives a constant remainder $C_{K,\rho}\alpha^{K+1}$, uniform on the
disk but confined to it, whereas (ii) is valid on the whole sector,
with a bound that decays like $|z|^{-1}$ as $|z|\to\infty$ but
degenerates as $z$ approaches the forbidden ray (i.e. as
$\delta\to0$). On the overlap either bound may be used; part (iii) is
precisely the statement that, together, they cover every compact
subset of $\mathbb C\setminus[1,\infty)$.

The angle $\theta=3\pi/4$ in Lemma \ref{lem-ML-Hankel} and
Theorem \ref{thm-ML-expansion} is merely a convenient choice.
Any $\theta\in(\pi/2,\pi)$ would give exponential decay of
$e^\zeta$ on the two rays while avoiding the branch cut.
Changing $\theta$ only affects the admissible value
$\alpha_0(\delta)=\delta/(2\theta)$ and the constants in
the estimates.
\end{rem}

The cases $K=0$ and $K=1$ are the ones used most often in the sequel: in the rate expansion of Proposition \ref{prop-rate} below, in the operator-valued arguments of Section \ref{subsec-rates}, and in the formal density expansion of Section \ref{subsec-density}. We record them explicitly.

\begin{cor}
{\rm (i)} ($K=0$.) The convergence
\begin{equation*}
\lim_{\alpha\rightarrow0^+}E_\alpha(z)=\frac{1}{1-z}
\end{equation*}
holds uniformly on compact subsets of $\mathbb{C}\setminus[1,\infty)$, with rate $O(\alpha)$. More precisely, for every $\rho\in(0,1)$ and $\delta\in(0,\pi)$ there exist $C_\rho,C_\delta>0$ such that
$$
\left|E_\alpha(z)-\frac{1}{1-z}\right|\leq C_\rho\,\alpha\qquad(|z|\leq\rho,\ 0<\alpha\leq1),
$$
$$
\left|E_\alpha(z)-\frac{1}{1-z}\right|\leq\frac{C_\delta\,\alpha\,|z|}{(1+|z|)\,|1-z|}\qquad(|\arg z|\geq\delta,\ 0<\alpha\leq\alpha_0(\delta)).
$$

{\rm (ii)} ($K=1$.) As $\alpha\rightarrow0^+$,
\begin{equation*}
E_\alpha(z)=\frac{1}{1-z}+\alpha\gamma\,\frac{z}{(1-z)^2}+O(\alpha^2).
\end{equation*}
More precisely, for every $\rho\in(0,1)$ and $\delta\in(0,\pi)$ there exist $C_\rho,C_\delta>0$ such that
$$
\left|E_\alpha(z)-\frac{1}{1-z}-\alpha\gamma\,\frac{z}{(1-z)^2}\right|\leq C_\rho\,\alpha^2\qquad(|z|\leq\rho,\ 0<\alpha\leq1),
$$
$$
\left|E_\alpha(z)-\frac{1}{1-z}-\alpha\gamma\,\frac{z}{(1-z)^2}\right|\leq\frac{C_\delta\,\alpha^2\,|z|}{|1-z|^2}\qquad(|\arg z|\geq\delta,\ 0<\alpha\leq\alpha_0(\delta)).
$$
\end{cor}

\begin{rem}\label{rem-ML-scope}
The qualitative limit is classical (Remark \ref{rem-ML-summability}).
Moreover, the endpoint coefficient
$\left.\partial_\alpha E_\alpha(z)\right|_{\alpha=0}
=\gamma z/(1-z)^2$ for $|z|<1$ is already recorded by
Apelblat \cite[Eq.~(99)]{Apelblat2020}, and Van Mieghem
\cite[art.~14, Eq.~(41)]{VanMieghem2021} gives the same formal
all-orders coefficient expansion in terms of Fermi--Dirac integrals.
Parameter differentiation at positive orders is also studied in
\cite{RGM2025}. Theorem \ref{thm-ML-expansion} supplies rigorous
finite-order remainder estimates, including the uniform sector
bounds needed for the operator-valued applications.
\end{rem}

Specializing the expansion to the solution curve $z=-\xi t^\alpha$ gives a full expansion of the scalar solution and, in particular, its exact first-order rate:

\begin{prop}\label{prop-rate}
Fix $\xi\in(-1,\infty)$. For every integer $K\geq0$ (with the same convention as in Theorem \ref{thm-ML-expansion}), the scalar solution curve admits the following expansion as $\alpha\rightarrow0^+$:
\begin{equation}\label{ML-curve-all-orders}
E_\alpha(-\xi t^\alpha)
=\frac{1}{1+\xi}+\sum_{k=1}^K\alpha^kP_k(\log t)\,\mathrm{Li}_{-k}(-\xi)
+O\Big(\alpha^{K+1}\bigl(1+|\log t|\bigr)^{K+1}\Big).
\end{equation}
The remainder is uniform for $t$ in every compact subset of $(0,\infty)$, with a constant depending on $\xi$, $K$ and the compact set. In particular, the case $K=1$ gives, for every $t>0$,
\begin{equation}\label{ML-rate}
E_\alpha(-\xi t^\alpha)=\frac{1}{1+\xi}-\alpha(\log t+\gamma)\,\frac{\xi}{(1+\xi)^2}+O(\alpha^2)\qquad(\alpha\rightarrow0^+),
\end{equation}
where the constant in $O(\alpha^2)$ depends on $\xi$ and $t$. For $\xi\neq0$ and $t\neq e^{-\gamma}$ the convergence in Lemma \ref{ML-limit} is of order exactly $\alpha$; for $\xi=0$ the error vanishes identically, and at the distinguished time $t=e^{-\gamma}$ the first-order term vanishes and the rate is $O(\alpha^2)$ (cf.\ Remark \ref{rem-expansion}(2)).
\end{prop}
\begin{proof}
The case $\xi=0$ is trivial, since $E_\alpha(0)=1$ and $\mathrm{Li}_{-k}(0)=0$. Let $\xi\neq0$ and apply Theorem \ref{thm-ML-expansion} to $z=-\xi t^\alpha$: for $\xi>0$ use part (ii), and for $\xi\in(-1,0)$ use part (i) after $\alpha$ is small enough. This gives
$$
E_\alpha(-\xi t^\alpha)
=\frac{1}{1+\xi t^\alpha}
+\sum_{j=1}^K\alpha^jc_j\,\mathrm{Li}_{-j}(-\xi t^\alpha)
+O(\alpha^{K+1}),
$$
with the required uniformity for $t$ in compact subsets of $(0,\infty)$.

Next, since $t^\alpha=e^{\alpha\log t}$, expand each $\mathrm{Li}_{-j}(-\xi t^\alpha)$ to order $K-j$ in $\alpha$, using
$$
\frac{d}{dv}\mathrm{Li}_{-j}(-\xi e^v)=\mathrm{Li}_{-(j+1)}(-\xi e^v).
$$
Taylor's theorem applies uniformly because $v$ between $0$ and $\alpha\log t$ stays in a compact range for $t$ in a compact interval and small $\alpha$. Substitution shows that the coefficient of $\alpha^k$ is
$$
\sum_{j=0}^k\frac{c_j(\log t)^{k-j}}{(k-j)!}\,\mathrm{Li}_{-k}(-\xi)
=P_k(\log t)\,\mathrm{Li}_{-k}(-\xi),
$$
where the constant term $1/(1+\xi t^\alpha)=1+\mathrm{Li}_0(-\xi t^\alpha)$ is included as the case $j=0$. The scalar remainder and the Taylor remainders together are $O(\alpha^{K+1}(1+|\log t|)^{K+1})$, which proves (\ref{ML-curve-all-orders}). For $K=1$, $P_1(\log t)=\log t+\gamma$ and $\mathrm{Li}_{-1}(-\xi)=-\xi/(1+\xi)^2$, giving (\ref{ML-rate}).
\end{proof}

%The first-order correction in (\ref{ML-rate}) is proportional to $\log t+\gamma$;
%its $t$-dependence shows once more that the convergence is not uniform near $t=0$, and identifies the boundary layer correction.

The scalar limit also yields the following vector-valued convergence of the subordination measures, which will be used both for abstract equations and for their stochastic interpretation.

\begin{lem}\label{lem-measure-weak}
For every bounded continuous $f:[0,\infty)\to X$,
\begin{equation}\label{measure-weak}
\int_0^\infty f(s)p_\alpha(s,t)\,ds
\longrightarrow\int_0^\infty f(s)e^{-s}\,ds
\end{equation}
in norm, uniformly for $t$ in compact subsets of $(0,\infty)$.
\end{lem}
\begin{proof}
Fix $[t_0,t_1]\subset(0,\infty)$. By (\ref{ML-Laplace}) and Lemma \ref{ML-limit} at $t=1$ (the local uniformity in $\xi$ absorbing the moving point $jt^\alpha$, which stays in a compact subset of $(-1,\infty)$ for $t\in[t_0,t_1]$ and small $\alpha$), for every integer $j\geq1$,
$$
\int_0^\infty e^{-js}p_\alpha(s,t)\,ds
=E_\alpha(-jt^\alpha)\longrightarrow\frac{1}{1+j}
=\int_0^\infty e^{-js}e^{-s}\,ds
$$
uniformly for $t\in[t_0,t_1]$. Taking finite linear combinations with vector coefficients gives the same convergence for every function $\sum_{j=1}^N e^{-js}x_j$, $x_j\in X$.

We next show that these finite sums are dense in $C_0([0,\infty);X)$. This density assertion is likely known, but we have not found a reference for it, so we include a detailed proof. Set $r=e^{-s}$, which identifies $C_0([0,\infty);X)$ with the space of functions $F\in C([0,1];X)$ satisfying $F(0)=0$. For such $F$, consider the Bernstein polynomials
$$
B_nF(r):=\sum_{k=0}^nF\Bigl(\frac{k}{n}\Bigr)\binom{n}{k}r^k(1-r)^{n-k}.
$$
Let $\omega_F(\delta):=\sup_{|u-v|\leq\delta}\|F(u)-F(v)\|$ be the modulus of continuity of $F$. Splitting the sum defining $B_nF(r)-F(r)$ according to $|k/n-r|<\delta$ or $|k/n-r|\geq\delta$, and using the identities $\sum_{k=0}^n\binom{n}{k}r^k(1-r)^{n-k}=1$ and $\sum_{k=0}^n(k/n-r)^2\binom{n}{k}r^k(1-r)^{n-k}=r(1-r)/n$, gives
$$
\|B_nF(r)-F(r)\|\leq\omega_F(\delta)+\frac{2\|F\|_\infty}{n\delta^2}.
$$
Letting first $n\to\infty$ and then $\delta\to0$ shows that $B_nF\to F$ uniformly on $[0,1]$. Since $B_nF(0)=F(0)=0$, each $B_nF(r)$ can be written as a finite sum $\sum_{j=1}^N r^jx_j$, as claimed.

For a general $f\in C_0([0,\infty);X)$, choose such a finite sum $g$ with $\|f-g\|_\infty$ arbitrarily small. Since $p_\alpha(s,t)\,ds$ and $e^{-s}\,ds$ are probability measures, replacing $f$ by $g$ changes either integral by at most $\|f-g\|_\infty$ in norm. The convergence already proved for $g$ therefore yields (\ref{measure-weak}) in norm, uniformly for $t\in[t_0,t_1]$, for every $f\in C_0([0,\infty);X)$.

The same assertion holds for every bounded continuous $f$. Indeed, differentiating (\ref{ML-Laplace}) at $\xi>0$ and letting $\xi\downarrow0$ (monotone convergence) gives
$$
\int_0^\infty s\,p_\alpha(s,t)\,ds
=\frac{t^\alpha}{\Gamma(1+\alpha)},
$$
which is uniformly bounded on $[t_0,t_1]$ for small $\alpha$. Let $f_R$ be the product of $f$ with a continuous cutoff that equals $1$ on $[0,R]$ and vanishes outside $[0,R+1]$. Applying (\ref{measure-weak}) to $f_R\in C_0([0,\infty);X)$, the two errors incurred by replacing $f_R$ by $f$ are bounded by
$$
\|f\|_\infty\left(\frac{t^\alpha}{R\Gamma(1+\alpha)}+e^{-R}\right),
$$
uniformly in $t$. Letting first $\alpha\to0^+$ and then $R\to\infty$ proves the assertion for $f$.
\end{proof}

\section{The limit as \texorpdfstring{$\alpha\rightarrow0^+$}{alpha to 0+}: the Caputo derivative and abstract fractional equations}\label{sec-abstract}

This section develops the limit theory for the abstract problem (\ref{FCP}). We first determine the limit of the Caputo derivative itself, then prove the convergence $S_\alpha(t)\phi\to(1-A)^{-1}\phi$ of the solutions, and finally derive operator-norm convergence rates and expansions, including an inner expansion resolving the boundary layer at $t=0$.

\subsection{The limit of the Caputo derivative as \texorpdfstring{$\alpha\rightarrow 0^+$}{alpha to 0+}}\label{sec-Caputo-limit}

As a first step we determine the limit of the fractional derivatives themselves; this is the point at which the Caputo and Riemann--Liouville derivatives part ways, and it dictates the form of the limit equation (\ref{limit-eq-intro}).

\begin{prop}\label{Caputo-limit-thm}
Let $T>0$.

(1) (Caputo derivative.) If $\phi\in C^1([0,T];X)$, then
\begin{equation}\label{limit-Caputo}
\lim_{\alpha\rightarrow 0^+}D_t^\alpha\phi(t)=\phi(t)-\phi(0),
\end{equation}
and the convergence is uniform on $[0,T]$. If $\phi\in W^{1,p}((0,T);X)$ for some $1\leq p<\infty$, then $D_t^\alpha\phi\rightarrow \phi-\phi(0)$ in $L^p((0,T);X)$ as $\alpha\rightarrow0^+$.

(2) (Riemann--Liouville derivative.) For $\phi$ as in {\rm(1)}, the Riemann--Liouville derivative satisfies instead
\begin{equation*}
\lim_{\alpha\rightarrow 0^+}{}_{RL}D_t^\alpha\phi(t)=\phi(t),
\end{equation*}
with the same modes of convergence, but for $\phi(0)\neq0$ the uniform convergence holds only on compact subsets of $(0,T]$, since $g_{1-\alpha}(t)\phi(0)$ is singular at $t=0$.
\end{prop}
\begin{proof}
We first record that
\begin{equation}\label{g1alpha-L1}
g_{1-\alpha}\longrightarrow 1\qquad\text{in }L^1(0,T)\ \text{as }\alpha\rightarrow0^+.
\end{equation}
Indeed, $g_{1-\alpha}(s)=s^{-\alpha}/\Gamma(1-\alpha)\rightarrow1$ for every $s>0$. For $0<\alpha\leq1/2$ one has $\Gamma(1-\alpha)\geq\Gamma(1)=1$, since $\Gamma$ is decreasing on $(0,1]$, and $s^{-\alpha}\leq 1+s^{-1/2}$. Hence $|g_{1-\alpha}(s)-1|\leq 2+s^{-1/2}\in L^1(0,T)$, and the dominated convergence theorem gives (\ref{g1alpha-L1}).

(1) By the definition of the Caputo derivative (Definition \ref{def-Caputo} with $m=1$) and $\int_0^t\phi'(s)ds=\phi(t)-\phi(0)$,
\begin{equation}\label{Caputo-error}
D_t^\alpha\phi(t)-(\phi(t)-\phi(0))=\int_0^t\left(g_{1-\alpha}(t-s)-1\right)\phi'(s)\,ds.
\end{equation}
Setting $M=\sup_{[0,T]}\|\phi'\|$, for every $t\in(0,T]$ this gives
$$
\|D_t^\alpha\phi(t)-(\phi(t)-\phi(0))\|\leq M\int_0^t|g_{1-\alpha}(s)-1|\,ds\leq M\,\|g_{1-\alpha}-1\|_{L^1(0,T)},
$$
which is independent of $t$ and tends to $0$ by (\ref{g1alpha-L1}). At $t=0$ both sides of (\ref{limit-Caputo}) vanish, by the convention of Definition \ref{def-Caputo}: $\|D_t^\alpha\phi(t)\|\leq Mt^{1-\alpha}/\Gamma(2-\alpha)\rightarrow0$ as $t\rightarrow0^+$. Hence the convergence is uniform on $[0,T]$. The $L^p$ convergence follows from (\ref{Caputo-error}) and Young's inequality on $(0,T)$:
$$
\|D_t^\alpha\phi-(\phi-\phi(0))\|_{L^p(0,T;X)}\leq \|g_{1-\alpha}-1\|_{L^1(0,T)}\,\|\phi'\|_{L^p(0,T;X)},
$$
which tends to $0$ by (\ref{g1alpha-L1}).

(2) Use the relation ${}_{RL}D_t^\alpha\phi(t)=D_t^\alpha\phi(t)+g_{1-\alpha}(t)\phi(0)$, part (1), and $g_{1-\alpha}(t)=t^{-\alpha}/\Gamma(1-\alpha)\rightarrow1$ uniformly on compact subsets of $(0,\infty)$ and in $L^p(0,T)$ for every $1\leq p<\infty$ (for the $L^p$ convergence dominate by $1+t^{-1/(2p)}\in L^p(0,T)$, valid for $\alpha\leq1/(2p)$).
\end{proof}

Thus, as $\alpha\rightarrow0^+$, the Caputo derivative converges to $\phi-\phi(0)$, whereas the Riemann--Liouville derivative converges to $\phi$. The term $-\phi(0)$ survives only in the Caputo limit and becomes a source term in the limit equation. For $\alpha\rightarrow1^-$, both derivatives converge to $\phi'$ locally uniformly for $t>0$.

\begin{example}\label{ex-exp}
Let $\phi(t)=e^{\lambda t}$. The series expansion of the exponential gives
$$
D_t^\alpha\phi(t)=t^{-\alpha}\sum_{m=1}^\infty\frac{(\lambda t)^m}{\Gamma(m+1-\alpha)}
\longrightarrow e^{\lambda t}-1=\phi(t)-\phi(0),
$$
the passage to the limit being justified as in Lemma \ref{ML-limit} (the same denominator bound applies since $m+1-\alpha\geq1$), while
$$
{}_{RL}D_t^\alpha\phi(t)=D_t^\alpha\phi(t)+\frac{t^{-\alpha}}{\Gamma(1-\alpha)}\,\phi(0)
\longrightarrow e^{\lambda t}=\phi(t)
$$
for every $t>0$. It is evident that the two limits differ exactly by the initial value $\phi(0)=1$. The extreme case is the constant function $\phi\equiv c$: its Caputo derivative vanishes identically, whereas ${}_{RL}D_t^\alpha\phi(t)=ct^{-\alpha}/\Gamma(1-\alpha)\rightarrow c$ for $t>0$. If $\phi(0)=0$ (for instance $\phi(t)=t$, where $D_t^\alpha\phi(t)={}_{RL}D_t^\alpha\phi(t)=t^{1-\alpha}/\Gamma(2-\alpha)\rightarrow t$), the two derivatives coincide and no difference appears.
\end{example}

\begin{rem}\label{rem-literature}
The limiting densities obtained in the literature for iterated Brownian motions fit naturally into the framework above, and it is worth clarifying how the formal step used in part of that literature relates to the exact limit equation (\ref{limit-eq-intro}). In \cite[Remark 3.4]{Orsingher2009}, the density of the $m$-times iterated Brownian motion converges to $u(x)=e^{-2|x|}$. It solves an equation of the form (\ref{FCP}) with order $2^{-m}$ and $A=2^{\,2^{-m}-2}\partial_x^2$. The same limiting density appears in the related models of \cite{Beghin2009}. To identify the limit equation, these papers formally set $\lim_{\alpha\rightarrow0^+}D_t^\alpha u=u$, an assumption stated explicitly there (see also \cite[Remark 4.4]{Beghin2009}). In \cite[Theorem 3.3]{Orsingher2017}, the limit is instead obtained and verified through the Fourier--Laplace transform, consistently with the present equation.
While the replacement $D_t^\alpha u\to u$ is legitimate for the Riemann--Liouville derivative (Proposition \ref{Caputo-limit-thm}(2)), the Caputo derivative used in these equations has the different limit $\phi-\phi(0)$. Nevertheless, for the delta initial data considered in those papers, the formally obtained limiting density is correct. Indeed, $u(x)=e^{-2|x|}$ satisfies, in the sense of distributions,
$$
u-\delta_0=\frac14\,\frac{\partial^2 u}{\partial x^2},
$$
because $\partial_x^2 e^{-2|x|}=4e^{-2|x|}-4\delta_0$. This is exactly our limit equation (\ref{limit-eq-intro}) with $A=\frac14\partial_x^2$ and $u(\cdot,0)=\delta_0$, where $\frac14=\lim_{m\rightarrow\infty}2^{\,2^{-m}-2}$ and $\delta_0$ accounts for the correction term $-u(\cdot,0)$. For regular initial data, however, the two procedures differ essentially: the equation $u=Au$ forces $u$ into the kernel of $I-A$, whereas the correct limit is $(1-A)^{-1}\phi$ (Theorem \ref{thm-homogeneous}).
\end{rem}

\subsection{Abstract equations: convergence of solutions}

With the notions of Section \ref{subsec-resolvent} in force, we now pass to the limit $\alpha\rightarrow0^+$ in the abstract equation.

\begin{thm}\label{thm-homogeneous}
Suppose $A$ generates a $\beta$-times resolvent family $S_\beta$ for some fixed $\beta\in(0,2]$ with $\sup_{t\geq0}\|S_\beta(t)\|\leq M<\infty$ for some $M>0$, and let $S_\alpha(t)$, $0<\alpha<\beta$, be the $\alpha$-times resolvent families given by the subordination formula (\ref{subordination-general}). Let $u_\alpha(t)=S_\alpha(t)\phi$ be the mild solution of
\begin{equation}\label{homogeneous-abstract}
\begin{cases}
D_t^\alpha u_\alpha(t)=Au_\alpha(t),  \quad t>0, \\
u_\alpha(0)=\phi.
\end{cases}
\end{equation}
Then for every $\phi\in X$ and every $t>0$,
\begin{equation}\label{homogeneous-limit}
\lim_{\alpha\rightarrow0^+}u_\alpha(t)=\lim_{\alpha\rightarrow0^+}S_\alpha(t)\phi=(1-A)^{-1}\phi=\int_0^\infty e^{-s}S_\beta(s)\phi\,ds,
\end{equation}
and the convergence is uniform on compact subsets of $(0,\infty)$. In other words, the solution of (\ref{homogeneous-abstract}) converges to the (constant in time) solution of the limit equation $u(t)-\phi=Au(t)$.
\end{thm}
\begin{proof}
For each $\phi\in X$, the function $s\mapsto S_\beta(s)\phi$ is bounded and continuous. Hence the subordination formula (\ref{subordination-general}) and Lemma \ref{lem-measure-weak}, applied with $\alpha$ replaced by $\alpha/\beta$, give
$$
u_\alpha(t)=S_\alpha(t)\phi
=\int_0^\infty p_{\alpha/\beta}(s,t)S_\beta(s)\phi\,ds
\longrightarrow\int_0^\infty e^{-s}S_\beta(s)\phi\,ds
$$
in norm, uniformly for $t$ in compact subsets of $(0,\infty)$. The Laplace identity (\ref{LT-resolventfamily}) at $\lambda=1$ gives $1\in\rho(A)$ and identifies the right-hand side as $(1-A)^{-1}\phi$. This proves (\ref{homogeneous-limit}).
\end{proof}

\begin{rem}\label{rem-boundary-layer}
(1) Since $S_\alpha(0)\phi=\phi$, convergence to $(1-A)^{-1}\phi$ fails at $t=0$, and cannot be uniform on an interval containing $0$, unless $\phi\in\ker A$. For $\phi\in\ker A$, $S_\alpha(t)\phi=\phi$ identically. The integrated families $\Sigma_\alpha(t)\phi:=\int_0^t S_\alpha(s)\phi\,ds$, however, converge to $t(1-A)^{-1}\phi$ uniformly on compact subsets of $[0,\infty)$, by Theorem \ref{thm-homogeneous} and the uniform bound; equivalently, this follows from the approximation theorem for vector-valued Laplace transforms \cite[Theorem 1.7.5, p.~41]{ABHN2011}.

(2) Theorem \ref{thm-homogeneous} also remains valid, with the same proof, when $A$ generates a $\beta$-times resolvent family $S_\beta$ for some $\beta\in(0,2]$ with $\|S_\beta(s)\|\leq Me^{\omega s}$ and $\omega<1$. The moment identity
$$
\int_0^\infty s^n\,p_\alpha(s,t)\,ds=\frac{n!\,t^{\alpha n}}{\Gamma(1+\alpha n)}
$$
and termwise integration of the positive series $e^{\omega' s}=\sum_n(\omega' s)^n/n!$ give $\int_0^\infty e^{\omega' s}p_\alpha(s,t)\,ds=E_\alpha(\omega't^\alpha)$ for every $\omega'>0$ (and hence, with (\ref{ML-Laplace}), for every real $\omega'$). Choose $\max\{0,\omega\}<\omega'<1$. Then
$$
\int_R^\infty e^{\omega s}p_\alpha(s,t)\,ds
\leq e^{-(\omega'-\omega)R}E_\alpha(\omega't^\alpha),
$$
whose right-hand side tends to zero as $R\to\infty$, uniformly for small $\alpha$ and $t$ in compact subsets of $(0,\infty)$. The cutoff argument in the proof of Lemma \ref{lem-measure-weak} therefore applies, and the limiting integral is $(1-A)^{-1}\phi$. The threshold cannot in general be extended to $\omega=1$, since $T(s)=e^sI$ gives the divergent scalar example $a=1$ below.

\end{rem}

To motivate the boundedness assumption, consider the scalar operator $A=aI$ on $\mathbb C$, with $a\in\mathbb R$. If $a<1$, then $E_\alpha(at^\alpha)\to(1-a)^{-1}$ by Lemma \ref{ML-limit} with $\xi=-a$. If $a\geq1$, positivity gives
$$
\liminf_{\alpha\to0^+}E_\alpha(at^\alpha)
\geq\sum_{n=0}^N a^n\qquad(N\geq0),
$$
so the solution diverges for every $t>0$. Thus $1\in\rho(A)$, which holds for $a>1$, is insufficient; likewise, exponential bounds with constants depending on $\alpha$ do not ensure the limit. Boundedness is a sufficient condition, not a necessary one, as $0<a<1$ also shows.

On the other hand, a single bounded resolvent family already yields convergence at the Laplace level. Suppose $A$ generates an $\alpha_0$-times resolvent family $S_{\alpha_0}$ that is bounded in operator norm. By subordination, the families for $0<\alpha<\alpha_0$ exist and are uniformly bounded by the same constant. Once $1\in\rho(A)$ is established,
$$
\int_0^\infty e^{-\lambda t}S_\alpha(t)\phi\,dt
=\lambda^{\alpha-1}(\lambda^\alpha-A)^{-1}\phi
\longrightarrow\lambda^{-1}(1-A)^{-1}\phi.
$$
The span of $\{e^{-\lambda t}:\lambda>0\}$ is dense in $L^1(0,\infty)$ \cite[Lemma 1.7.1, p.~40]{ABHN2011}. Together with the uniform bound, this extends the convergence, in norm, to every scalar test function in $L^1$. The condition $1\in\rho(A)$ and the associated positive-axis resolvent estimates are already part of the generation theory \cite[Section 2.1]{Baj}. Theorem \ref{thm-unifbdd} later upgrades this Laplace-level convergence to operator-norm convergence with rate $O(\alpha)$ by establishing sectoriality.

\subsection{Operator-norm convergence rates and expansions}\label{subsec-rates}

We now turn to the convergence rates. The appropriate class of generators is then that of sectorial operators---and this is not an additional restriction: by Theorem \ref{thm-unifbdd} below, a single bounded fractional resolvent family alone already forces $A$ to be sectorial. For such operators the convergence takes place in operator norm, with rate $O(\alpha)$, and indeed admits a complete asymptotic expansion in powers of $\alpha$. The proofs below use the quantitative expansion estimates for the Mittag-Leffler function established in Section \ref{subsec-ML-expansion}.

We first fix the precise functional-calculus setting. An operator $A$ is called sectorial of angle $\theta\in(0,\pi)$ around the negative real axis if $\sigma(A)\subset\overline{\Sigma_\theta}$, where $\Sigma_\theta=\{z\in\mathbb{C}\setminus\{0\}:|\arg(-z)|<\theta\}$, and
\begin{equation}\label{sectorial-est}
\|(\zeta-A)^{-1}\|\leq\frac{M}{|\zeta|},\qquad \zeta\notin\overline{\Sigma_\theta},\ \zeta\neq0.
\end{equation}
In particular, $1\in\rho(A)$ and $\|(1-A)^{-1}\|\leq M$, since $1\notin\overline{\Sigma_\theta}$. For the definition of sectorial operators, we refer to \cite[Section 3.8, p.~162]{ABHN2011}; for the related fractional resolvent calculus, see \cite{LiChenLi2010}.

For the Dunford representation below, set $\varphi=(\pi-\theta)/2$ and let $\Gamma=\{re^{\pm i\varphi}:r\geq\tfrac12\}\cup\{\tfrac12 e^{i\psi}:|\psi|\leq\varphi\}$. The contour is traversed from the lower ray to the upper one. It winds once counterclockwise around the region $\{\zeta:|\zeta|\leq\tfrac12\}\cup\{\zeta:|\arg(-\zeta)|<\pi-\varphi\}$, which contains $\overline{\Sigma_\theta}$. The arc passes on the side of the origin {\em opposite} to the spectral sector. Moreover, every point of $\Gamma$ satisfies $|\arg(-\zeta)|\geq\pi-\varphi=\frac{\pi+\theta}{2}>\theta$, so that (\ref{sectorial-est}) applies along the whole contour.

For such operators, the $\alpha$-times resolvent family exists for $0<\alpha<2(1-\theta/\pi)$. Moreover, whenever $\pi\alpha/2<\varphi<\pi-\theta$ (which holds when $\alpha\leq\alpha_0(\varphi)=2\varphi/(3\pi)$), it admits the Dunford representation
\begin{equation}\label{sectorial-dunford}
S_\alpha(t)=\frac{1}{2\pi i}\int_\Gamma E_\alpha(\zeta t^\alpha)(\zeta-A)^{-1}\,d\zeta,\qquad t>0,
\end{equation}
with the integral converging in operator norm. This representation is known; see \cite[Remark 3.2(c)]{LiChenLi2010}, with the sign of the generator changed to our convention. It does not require $0\in\rho(A)$.

The shape of the expansion is predicted on the resolvent side: the formal expansion
$$
\lambda^{\alpha-1}(\lambda^\alpha-A)^{-1}=\lambda^{-1}(1-A)^{-1}+\alpha\log\lambda\,\bigl[\lambda^{-1}(1-A)^{-1}-\lambda^{-1}(1-A)^{-2}\bigr]+O(\alpha^2(\log\lambda)^2)
$$
suggests $\log t$-type corrections for $S_\alpha(t)\phi-(1-A)^{-1}\phi$, and the scalar case of Proposition \ref{prop-rate} already indicates that such a first-order term is sharp. To state the expansion, recall that $\mathrm{Li}_{-k}(z)=\frac{zQ_k(z)}{(1-z)^{k+1}}$ with $Q_k$ the Eulerian polynomials ($k\geq1$), and that the polynomials $P_k$ were introduced in Section \ref{subsec-ML-expansion}.

The proof has two steps. First, Theorem \ref{thm-ML-expansion} is applied to the scalar factor $E_\alpha(\zeta t^\alpha)$ on $\Gamma$. Second, each $\mathrm{Li}_{-j}(\zeta e^{\alpha\log t})$ is expanded in the variable $v=\alpha\log t$ and integrated termwise. The Dunford calculus identifies the coefficients, while (\ref{sectorial-est}) makes the common remainder integrable in operator norm.

\begin{thm}\label{thm-sectorial}
Let $A$ be sectorial of angle $\theta\in(0,\pi)$ around the negative real axis, and let $S_\alpha$ be the $\alpha$-times resolvent family in (\ref{sectorial-dunford}). Then for every integer $K\geq0$ (same convention as in Theorem \ref{thm-ML-expansion}), the following expansion holds as $\alpha\rightarrow0^+$, in operator norm and uniformly on compact subsets of $(0,\infty)$,
\begin{equation}\label{sol-higher-order}
S_\alpha(t)=(1-A)^{-1}+\sum_{k=1}^K\alpha^k\,P_k(\log t)\,\mathrm{Li}_{-k}(A)+O\Big(\alpha^{K+1}\bigl(1+|\log t|\bigr)^{K+1}\Big),
\end{equation}
where $\mathrm{Li}_{-k}(A)=A\,Q_k(A)\,(1-A)^{-k-1}$. The implied constant and the admissible range of $\alpha$ depend only on $M$, $\theta$, $K$ and the compact set. The estimate (\ref{sol-higher-order}) also holds uniformly whenever $|\alpha\log t|\leq1/4$, with a constant depending only on $M$, $\theta$ and $K$.
\end{thm}
When $A=-\xi$ on $\mathbb C$, (\ref{sol-higher-order}) reduces to the scalar solution-curve expansion (\ref{ML-curve-all-orders}) in Proposition \ref{prop-rate}.
\begin{proof}
On the rays, Theorem \ref{thm-ML-expansion}(ii) gives $E_\alpha(\zeta t^\alpha)=O(|\zeta|^{-1})$ for fixed $t>0$; together with (\ref{sectorial-est}), this proves norm convergence of the contour integral.

Fix $[t_0,t_1]\subset(0,\infty)$ and $K\geq0$. Choose $\alpha_1\leq\alpha_0(\varphi)$ so that $\max\{1,t_1^{\alpha_1}\}\leq3/2$. Apply Theorem \ref{thm-ML-expansion} to $\zeta t^\alpha$: on the arc use part (i) with $\rho=3/4$, and on the rays use part (ii). Next expand each $\mathrm{Li}_{-j}(\zeta e^{\alpha\log t})$, $0\leq j\leq K$, to order $K-j$, using
$$
\frac{d}{dv}\mathrm{Li}_{-j}(\zeta e^v)
=\mathrm{Li}_{-(j+1)}(\zeta e^v).
$$
For $v$ between $0$ and $\alpha\log t$, $e^v$ stays between $\min\{1,t_0^{\alpha_1}\}$ and $\max\{1,t_1^{\alpha_1}\}$. On the rays,
$$
|1-\zeta e^v|\geq\sin(\varphi/2)(1+|\zeta|e^v),
$$
while on the arc $|\zeta e^v|\leq3/4$. The Taylor remainders in $v$ are handled by the same two bounds: since $\frac{d}{dv}\mathrm{Li}_{-j}(\zeta e^v)=\mathrm{Li}_{-(j+1)}(\zeta e^v)$, and the rational form of $\mathrm{Li}_{-m}$ gives $|\mathrm{Li}_{-m}(w)|\leq C_m|w|/(1+|w|)^2$ whenever $|1-w|\geq\sin(\varphi/2)(1+|w|)$ or $|w|\leq3/4$, while $e^v$ stays bounded above and below, every derivative is $O\bigl(|\zeta|/(1+|\zeta|)^2\bigr)$ uniformly in $v$. Thus the scalar remainder and all Taylor remainders have the common bound
$$
C\alpha^{K+1}(1+|\log t|)^{K+1}
\frac{|\zeta|}{(1+|\zeta|)^2},\qquad \zeta\in\Gamma.
$$
The constant term is $1+\mathrm{Li}_0(\zeta)=(1-\zeta)^{-1}$, and the coefficient of $\alpha^k$ is
$$
\sum_{j=0}^k\frac{c_j(\log t)^{k-j}}{(k-j)!}\,
\mathrm{Li}_{-k}(\zeta)
=P_k(\log t)\,\mathrm{Li}_{-k}(\zeta).
$$

These rational functions decay at infinity and have their only pole at $1$, outside the region enclosed by $\Gamma$. The Dunford calculus therefore identifies their contour integrals with $(1-A)^{-1}$ and $\mathrm{Li}_{-k}(A)$. The remainder is integrable in operator norm because
$$
\int_\Gamma \frac{|\zeta|}{(1+|\zeta|)^2}
\|(\zeta-A)^{-1}\|\,|d\zeta|
\leq M\int_\Gamma\frac{|d\zeta|}{(1+|\zeta|)^2}<\infty.
$$
This proves (\ref{sol-higher-order}). The same proof gives a constant depending only on $M,\theta,K$ whenever $|\alpha\log t|\leq1/4$: then $e^v\in[e^{-1/4},e^{1/4}]$ and the arc still lies in $|\zeta e^v|<3/4$. This additional uniformity will justify the overlap calculation in Remark \ref{rem-inner-matching}.
\end{proof}

\begin{cor}\label{cor-rate-expansion}
Under the hypotheses of Theorem \ref{thm-sectorial}:

{\rm (i)} $\|S_\alpha(t)-(1-A)^{-1}\|=O(\alpha)$ in operator norm, uniformly on compact subsets of $(0,\infty)$;

{\rm (ii)} the first-order expansion reads, as $\alpha\rightarrow0^+$, in operator norm and uniformly on compact subsets of $(0,\infty)$,
\begin{equation}\label{sol-expansion}
S_\alpha(t)=(1-A)^{-1}+\alpha(\log t+\gamma)\,A(1-A)^{-2}+O\bigl(\alpha^2(1+(\log t)^2)\bigr),
\end{equation}
where $\gamma$ is the Euler constant. By Theorem \ref{thm-unifbdd} below, both conclusions hold whenever $A$ generates a bounded fractional resolvent family of some order.
\end{cor}
\begin{proof}
Part (i) is the case $K=0$ of (\ref{sol-higher-order}), the factor $1+|\log t|$ being bounded on compact subsets of $(0,\infty)$; part (ii) is the case $K=1$, since $P_1(\log t)=\log t+\gamma$ and $\mathrm{Li}_{-1}(A)=A(1-A)^{-2}$.
\end{proof}

\begin{rem}
For $\theta<\pi/2$ the operator $A$ generates a bounded analytic semigroup, and strong convergence $S_\alpha(t)\phi\rightarrow(1-A)^{-1}\phi$ already follows from Theorem \ref{thm-homogeneous}; Corollary \ref{cor-rate-expansion}(i) upgrades this to norm convergence with the rate $O(\alpha)$. The range $\theta\in[\pi/2,\pi)$ is genuinely beyond the semigroup setting: let $X=L^2(\mathbb{R})$ and $(Af)(\omega)=m(\omega)f(\omega)$ with $m(\omega)=(1+|\omega|)\,e^{i\frac{\pi}{3}\,\mathrm{sign}\,\omega}$ and maximal domain $D(A)=\{f\in L^2(\mathbb R):mf\in L^2(\mathbb R)\}$. Then $\sigma(A)$ is the union of the two rays $\{re^{\pm i\pi/3}:r\geq1\}$; hence $A$ is sectorial of any angle $\theta>2\pi/3$ with the estimate (\ref{sectorial-est}), $1\in\rho(A)$, but $\sup\mathrm{Re}\,\sigma(A)=+\infty$, so $A$ generates no $C_0$-semigroup. Corollary \ref{cor-rate-expansion}(i) still applies and yields $\|S_\alpha(t)-(1-A)^{-1}\|=O(\alpha)$ for all sufficiently small $\alpha>0$. Moreover, the representation (\ref{sectorial-dunford}) holds whenever $\alpha<2(1-\theta/\pi)$. In this example, choosing $\theta$ arbitrarily close to $2\pi/3$ from above covers every $\alpha<2/3$; see the angle condition $\pi\alpha/2<\omega<\pi-\theta$ in \cite[Remark 3.2(c)]{LiChenLi2010}.
\end{rem}

As announced, the sectoriality assumption of Theorem \ref{thm-sectorial} is in fact a consequence of a single bounded fractional resolvent family. The positive-axis resolvent estimates are classical \cite[Section 2.1]{Baj}; the point here is to upgrade them to an explicit sectorial estimate and thereby answer in the affirmative the Laplace-transform problem raised in the discussion above.

\begin{thm}\label{thm-unifbdd}
Assume that $A$ generates, for some $\alpha_0>0$, an $\alpha_0$-times resolvent family $S_{\alpha_0}$ with $\sup_{t\geq0}\|S_{\alpha_0}(t)\|\leq M<\infty$. Then $A$ is sectorial around the negative real axis of angle $\pi-\theta_0$ with $\theta_0=2\arcsin\frac{1}{2M}>0$, and the conclusions of Theorem \ref{thm-sectorial} hold: in particular, for every compact $[t_0,t_1]\subset(0,\infty)$ there exist $\alpha_1>0$, depending only on $\alpha_0$, $M$ and $t_1$, and $C>0$, depending only on $M$, $t_0$ and $t_1$, such that
$$
\|S_\alpha(t)-(1-A)^{-1}\|\leq C\,\alpha,\qquad t\in[t_0,t_1],\ 0<\alpha\leq\alpha_1.
$$
\end{thm}
\begin{proof}
By \cite[Section 2.1, (2.5)--(2.6) and Theorem 2.8, (2.12)]{Baj}, the boundedness assumption gives the positive-axis resolvent estimate
$$
(0,\infty)\subset\rho(A),\qquad
\|(\mu-A)^{-1}\|\leq M/\mu\qquad(\mu>0).
$$
For $0<\alpha<\alpha_0$, the subordination formula (\ref{subordination-general}) supplies the families $S_\alpha$ and yields
\begin{equation}\label{unifbdd-hyp}
\sup_{\alpha\in(0,\alpha_0],\,t\geq0}\|S_\alpha(t)\|\leq M.
\end{equation}
Also $M\geq1$, since $S_\alpha(0)=I$, and hence
$$
0<\theta_0=2\arcsin\frac{1}{2M}\leq\frac{\pi}{3},
\qquad M\sin\theta_0=\sqrt{1-\frac{1}{4M^2}}<1.
$$
For $\zeta\neq0$ with $|\arg\zeta|\leq\theta_0$, choose
$\mu=|\zeta|/\cos(\arg\zeta)>0$. Then
$$
\frac{|\zeta-\mu|}{\mu}=|\sin(\arg\zeta)|\leq\sin\theta_0,
$$
so the Neumann series gives $\zeta\in\rho(A)$ and
$$
\|(\zeta-A)^{-1}\|
\leq\frac{M}{\mu(1-M\sin\theta_0)}
\leq\frac{M}{|\zeta|(1-M\sin\theta_0)}.
$$
Thus (\ref{sectorial-est}) holds with angle $\pi-\theta_0$ and a constant depending only on $M$. That a single bounded fractional resolvent family forces sectoriality of $A$ is a known structural fact (cf.\ \cite[Remark 2.8(b)]{LiChenLi2010}, where the angle may depend on $\alpha_0$); the Neumann series argument above yields an angle depending only on $M$. Apply Theorem \ref{thm-sectorial}; its families agree with the given $S_\alpha$, which satisfy (\ref{unifbdd-hyp}), by uniqueness of the Laplace transform. The proof of that theorem permits $\alpha_1$ to depend only on the sector angle and $t_1$, with the remainder constant also depending on $t_0$; restricting $\alpha_1\leq\alpha_0$ yields the asserted dependence.
\end{proof}

When $\alpha_0=1$, the angle obtained above satisfies $\pi-\theta_0\geq2\pi/3$. In our convention this is much weaker than the sectorial angle $<\pi/2$ associated with bounded analytic semigroups, so it does not assert bounded analyticity and is consistent with classical semigroup theory.

\begin{cor}\label{cor-semigroup-rate}
If $A$ generates a bounded $C_0$-semigroup $(T(t))_{t\geq0}$, then the resolvent families $S_\alpha$ given by (\ref{subordination2}) satisfy $\|S_\alpha(t)-(1-A)^{-1}\|=O(\alpha)$ in operator norm, uniformly on compact subsets of $(0,\infty)$. Indeed, the subordination formula (\ref{subordination2}) gives $\|S_\alpha(t)\|\leq\sup_{s\geq0}\|T(s)\|<\infty$, so Theorem \ref{thm-unifbdd} applies. Thus the strong convergence of Theorem \ref{thm-homogeneous} is in fact norm convergence with a linear rate, and the full expansion (\ref{sol-higher-order}) is available.
\end{cor}

\begin{rem}\label{rem-expansion}
(1)  For $A=aI$, $a<1$, (\ref{sol-expansion}) becomes
$$
E_\alpha(at^\alpha)=(1-a)^{-1}
+\alpha(\log t+\gamma)a(1-a)^{-2}+O(\alpha^2),
$$
as in Proposition \ref{prop-rate}, which covers every $a<1$. Thus both the scalar coefficient and the operator factor $A(1-A)^{-2}$ are sharp.

(2) The first-order term vanishes at $t_*=e^{-\gamma}\approx0.5615$, where the scalar factor $\log t+\gamma$ changes sign. The first-order term also vanishes on $\ker A$: the resolvent equation gives $S_\alpha(t)\phi=\phi$ for every $\phi\in\ker A$, so convergence is exact there.

(3) The Riemann--Liouville integrals form a semigroup in the order,
$$
I_{a+}^{\alpha}I_{a+}^{\beta}=I_{a+}^{\alpha+\beta},\qquad I_{a+}^{0}=I.
$$
For sufficiently regular $\varphi$ and $x>a$, its infinitesimal generator appears in the first-order expansion
$$
(I_{a+}^\alpha\varphi)(x)=\varphi(x)+\alpha\left[
\gamma\varphi(x)+\frac{d}{dx}\int_a^x\log(x-t)\varphi(t)\,dt\right]+o(\alpha).
$$
The expression in brackets is the generator of this order semigroup; see \cite[Remark 2.5]{SKM1993} and \cite[Section 23.16]{HP1957} for its identification and domain. The Euler constant in this generator and in (\ref{sol-expansion}) has the same source, namely $1/\Gamma(1+\alpha)=1+\gamma\alpha+O(\alpha^2)$: it enters through the normalization of the Riemann--Liouville kernel in one case and through the coefficients $1/\Gamma(1+\alpha n)$ of the Mittag--Leffler series in the other.

At the level of the fractional operators themselves, the appropriate comparison is between $I^\alpha$ and $D_t^\alpha$. The former tends to the identity as $\alpha\to0^+$, whereas the Caputo derivative tends to $\phi-\phi(0)$. This discrepancy is inherited by the solution families and accounts for their resolvent endpoint $(1-A)^{-1}$. The family $\{S_\alpha(t)\}$ is not itself a semigroup in the order, so the Hille--Phillips generator theory does not apply directly; the Laplace and Dunford methods used above serve as substitutes.

\end{rem}

The expansion (\ref{sol-higher-order}) degenerates as $t\rightarrow0^+$ through the factor $\log t$. The boundary layer is resolved by the inner variable $s=-\alpha\log t$, i.e.\ $t=e^{-s/\alpha}$, for which $t^\alpha=e^{-s}$ exactly.

\begin{prop}\label{thm-inner}
Under the hypotheses of Theorem \ref{thm-sectorial}, for every integer $K\geq0$ (same convention as in Theorem \ref{thm-ML-expansion}), the following inner expansion holds as $\alpha\rightarrow0^+$:
\begin{equation}\label{sol-inner}
S_\alpha\big(e^{-s/\alpha}\big)=(1-e^{-s}A)^{-1}+\sum_{k=1}^K\alpha^k\,c_k\,\mathrm{Li}_{-k}\big(e^{-s}A\big)+O(\alpha^{K+1}),
\end{equation}
in operator norm, uniformly for $s\in[0,\infty)$. The implied constant depends only on $M$, $\theta$ and $K$, and the admissible range of $\alpha$ only on $\theta$. In particular, the inner profile $\Phi(s)=(1-e^{-s}A)^{-1}$ interpolates between the outer limit $\Phi(0)=(1-A)^{-1}$ and the initial datum: $\Phi(s)\rightarrow I=S_\alpha(0)$ strongly as $s\rightarrow\infty$.
\end{prop}
\begin{proof}
Since $t^\alpha=e^{-s}$ exactly, (\ref{sectorial-dunford}) reads $S_\alpha(e^{-s/\alpha})=\frac{1}{2\pi i}\int_\Gamma E_\alpha(\zeta e^{-s})(\zeta-A)^{-1}\,d\zeta$. Apply Theorem \ref{thm-ML-expansion} with $z=\zeta e^{-s}$: on the arc by part (i), on the rays by part (ii) with $\delta=\varphi$, noting that $|1-\zeta e^{-s}|\geq\sin(\varphi/2)(1+|\zeta|e^{-s})$ (Remark \ref{rem-sector-equivalent}). The pole $e^s\geq1$ of $\zeta\mapsto\mathrm{Li}_{-k}(e^{-s}\zeta)$ lies outside the region enclosed by $\Gamma$, so the Dunford calculus identifies the coefficients as $c_k\,\mathrm{Li}_{-k}(e^{-s}A)$. The remainder is uniform in $s\geq0$: on the rays it is bounded by
$$
C\alpha^{K+1}\int_{1/2}^\infty\frac{re^{-s}}{(1+re^{-s})^2}\cdot\frac{M}{r}\,dr
\leq CM\alpha^{K+1}\int_0^\infty\frac{du}{(1+u)^2},
$$
with $u=re^{-s}$, and the arc contributes $O(\alpha^{K+1})$ directly. Moreover, $(1-e^{-s}A)^{-1}=e^s(e^s-A)^{-1}$ and $e^s>0$ lies outside $\overline{\Sigma_\theta}$, so (\ref{sectorial-est}) gives $\|(1-e^{-s}A)^{-1}\|\leq M$. Finally, for $\phi\in D(A)$,
$$
\|(1-e^{-s}A)^{-1}\phi-\phi\|
= e^{-s}\|(1-e^{-s}A)^{-1}A\phi\|
\leq M e^{-s}\|A\phi\|\to0.
$$
The density of $D(A)$ extends this convergence to every $\phi\in X$.
\end{proof}

\begin{rem}\label{rem-inner-matching}
(1)  Matching follows from finite Taylor expansions at $s=0$, using
$$
\frac{d}{ds}\mathrm{Li}_{-j}(e^{-s}z)
=-\mathrm{Li}_{-(j+1)}(e^{-s}z).
$$
The rational functional calculus and the estimates in the proof of Theorem \ref{thm-sectorial} give, for $0\leq s\leq1/4$,
\begin{align*}
&(1-e^{-s}A)^{-1}+\sum_{j=1}^K\alpha^jc_j\mathrm{Li}_{-j}(e^{-s}A)\\
&\quad=(1-A)^{-1}+\sum_{k=1}^K\alpha^kP_k(-s/\alpha)\mathrm{Li}_{-k}(A)
+O\left(\sum_{j=0}^K\alpha^j s^{K-j+1}\right).
\end{align*}
On the overlap scale $s=\alpha^\sigma$, $0<\sigma<1$, this last error is $O(\alpha^{(K+1)\sigma})$. The outer estimate is valid there by its uniformity for $|\alpha\log t|\leq1/4$, established in Theorem \ref{thm-sectorial}; its remainder has the same order, while the inner remainder is $O(\alpha^{K+1})$. This justifies matching without assuming convergence of either infinite asymptotic series.

(2)  For $t=t(\alpha)\to0$, the outer leading term remains valid when $\alpha|\log t|\to0$. On the scale $t=e^{-s/\alpha}$ with fixed $s>0$, the profile instead retains its $s$-dependence through $(1-e^{-s}A)^{-1}$. This is the exponential scale of the initial layer.

(3)  At $s=0$, (\ref{sol-inner}) agrees with (\ref{sol-higher-order}) at $t=1$, since $P_k(0)=c_k$. Every correction vanishes on $\ker A$. At the density level ($A=-\partial_x$, Section \ref{subsec-density}), the inner profile is $e^se^{-xe^s}$, an exponential density with mean $e^{-s}$; consistently, for $t=e^{-s/\alpha}$,
$$
\mathbb EH_\alpha(t)=\frac{e^{-s}}{\Gamma(1+\alpha)}\longrightarrow e^{-s}.
$$
\end{rem}

We conclude this section with the nonhomogeneous equation, whose treatment uses the operator-norm convergence of Corollary \ref{cor-semigroup-rate} in an essential way.

\begin{thm}\label{thm-nonhomogeneous}
Suppose $A$ generates a bounded $C_0$-semigroup $(T(t))_{t\geq0}$ on $X$, and let $f\in C([0,\infty);X)$. Then the mild solution of
$$
\begin{cases}
D_t^\alpha u_\alpha(t)=Au_\alpha(t)+f(t),  \quad t>0, \\
u_\alpha(0)=\phi
\end{cases}
$$
satisfies, for every $t>0$,
\begin{equation}\label{nonhomogeneous-limit}
\lim_{\alpha\rightarrow0^+}u_\alpha(t)=(1-A)^{-1}(\phi+f(t)),
\end{equation}
uniformly on compact subsets of $(0,\infty)$. In contrast with the homogeneous case, the limit can depend on $t$.
\end{thm}
\begin{proof}
By the Volterra variation-of-parameters formula \cite[Section 1.2, Proposition 1.2, p.~33]{Pruss}, the mild solution takes the form
$$
u_\alpha(t)=S_\alpha(t)\phi+\int_0^t h_\alpha(r)f(t-r)\,dr,
$$
where the forcing kernel $h_\alpha$ is determined by $\widehat h_\alpha(\lambda)=(\lambda^\alpha-A)^{-1}$. The variation-of-constants formula and this forcing kernel are known in the theory of fractional evolution equations. For the reader's convenience, we include a detailed derivation in the present subordination setting before analyzing the forcing term.

Write $M=\sup_{s\geq0}\|T(s)\|$ and introduce the forcing kernel, with the integral understood strongly,
$$
h_\alpha(r):=\frac{\alpha}{r}\int_0^\infty
s\,p_\alpha(s,r)T(s)\,ds,\qquad r>0.
$$
The first-moment formula gives
$$
\|h_\alpha(r)\|\leq\frac{M\alpha r^{\alpha-1}}{\Gamma(1+\alpha)}
=M g_\alpha(r).
$$
Moreover, (\ref{inverse-stable-pdf-Laplace}) and Fubini's theorem yield, for $s,\lambda>0$,
$$
\int_0^\infty e^{-\lambda r}\frac{\alpha s}{r}p_\alpha(s,r)\,dr
=\alpha s\int_\lambda^\infty \mu^{\alpha-1}e^{-s\mu^\alpha}\,d\mu
=e^{-s\lambda^\alpha}.
$$
Consequently $\widehat h_\alpha(\lambda)=(\lambda^\alpha-A)^{-1}$. Comparison with (\ref{LT-resolventfamily}) verifies the mild formula displayed above and gives the integrated kernel
$$
\int_0^t h_\alpha(r)\,dr=(g_\alpha*S_\alpha)(t).
$$
The formula for $u_\alpha$ is first verified for bounded $f$ by Laplace transforms and then for general continuous $f$ by a cutoff in time.

Fix $[t_0,t_1]\subset(0,\infty)$. Since
$\int_0^t g_\alpha(r)\,dr=t^\alpha/\Gamma(1+\alpha)\to1$
uniformly on this interval and $g_\alpha$ concentrates at $0$, Corollary \ref{cor-semigroup-rate} gives
$$
(g_\alpha*S_\alpha)(t)\longrightarrow(1-A)^{-1}
$$
in operator norm, uniformly on $[t_0,t_1]$. To see the uniformity directly, split the convolution at $r=t-\delta$, with $0<\delta<t_0/2$. On the first part $t-r\geq\delta$, so the operator-norm convergence is uniform; the remaining kernel mass is at most
$\int_{t_0/2}^{t_1}g_\alpha(r)\,dr\to0$.

Finally, write the forcing term as
$$
\int_0^t h_\alpha(r)f(t-r)\,dr
=(g_\alpha*S_\alpha)(t)f(t)
+\int_0^t h_\alpha(r)\bigl(f(t-r)-f(t)\bigr)\,dr.
$$
The last integral tends to zero uniformly on $[t_0,t_1]$: on $0<r<\delta$ use the uniform continuity of $f$ on $[0,t_1]$ and the uniformly bounded mass of $g_\alpha$, while on $r\geq\delta$ use $\|h_\alpha(r)\|\leq M g_\alpha(r)$ and $\int_\delta^{t_1}g_\alpha(r)\,dr\to0$. Together with Theorem \ref{thm-homogeneous}, this proves (\ref{nonhomogeneous-limit}).
\end{proof}

\section{Concrete operators and the stochastic interpretation}\label{sec-concrete}

This section specializes the abstract theory to four concrete right-hand sides: the operator $-\partial_x$ behind the density of the inverse stable subordinator, the Laplacian on the whole space and on bounded domains, the space-time fractional operator $-(-\Delta)^\beta$, and the difference operator governing the fractional Poisson process. The closing subsection (Section \ref{sec-potential}) develops the stochastic interpretation from a different angle: the limit solution is the $1$-potential of the underlying Markov process, related to the Feynman--Kac formula and to a half-space elliptic extension.

\subsection{The operator \texorpdfstring{$-\partial_x$}{-partial\_x}: the density of the inverse stable subordinator}\label{subsec-density}

We begin with the density $p_\alpha$ of the inverse stable subordinator itself, which corresponds to the choice $A=-\partial_x$ with a delta initial datum. We record its limit and then give a formal first-order expansion suggested by the Mittag--Leffler expansion.

\begin{prop}\label{H-convergence}
Let $H_\alpha(t)$ be the inverse $\alpha$-stable subordinator and let $\tau$ be a unit exponential random variable. Then $H_\alpha(t)\Rightarrow\tau$ (convergence in distribution) as $\alpha\rightarrow0^+$ for every $t>0$, uniformly for $t$ in compact subsets of $(0,\infty)$: for $\alpha_j\rightarrow0^+$ and $t_j\rightarrow t_0>0$ one has $H_{\alpha_j}(t_j)\Rightarrow\tau$. Equivalently, for every $t>0$ the density $p_\alpha(\cdot,t)$ converges weakly to the unit exponential density $e^{-x}$ on $[0,\infty)$.
\end{prop}
\begin{proof}
Since $p_\alpha(\cdot,t)$ is the density of $H_\alpha(t)$, the scalar case $X=\mathbb R$ of Lemma \ref{lem-measure-weak} gives precisely the weak convergence of these laws to the unit exponential law, locally uniformly in $t$. The sequential formulation in the proposition is equivalent to this local uniformity.
\end{proof}

\begin{rem}\label{rem-p-alpha}
(1) Formally, the limit can also be read off from the time Laplace transform (\ref{inverse-stable-pdf-Laplace}): $\lambda^{\alpha-1}e^{-x\lambda^\alpha}\rightarrow\lambda^{-1}e^{-x}=\int_0^\infty e^{-\lambda t}e^{-x}\,dt$.

(2) The limit density $e^{-x}$ is independent of $t$. At the level of equations this is explained as follows: $p_\alpha(x,t)$ solves the time fractional equation \cite{Meerschaert2012book}
\begin{equation*}
\begin{cases}
D_t^\alpha p_\alpha(x,t)=-\frac{\partial p_\alpha(x,t)}{\partial x},  \quad t>0,\ x>0,\\
p_\alpha(x,0)=\delta_0(x).
\end{cases}
\end{equation*}
The limit equation is $p(x,t)-\delta_0(x)=-\partial_x p(x,t)$, whose solution is $p(x,t)=e^{-x}$, independent of $t$; indeed $\partial_x(e^{-x}\mathbf{1}_{\{x>0\}})=-e^{-x}+\delta_0$ in the sense of distributions, so that $p-\delta_0=-\partial_x p$ holds.

(3) In contrast, the $\alpha$-stable subordinator $Y_\alpha$ degenerates as $\alpha\rightarrow0^+$: since $\mathbb{E}e^{-\lambda Y_\alpha(s)}=e^{-s\lambda^\alpha}\rightarrow e^{-s}$ for every $\lambda>0$, part of the mass escapes to infinity. Accordingly, the normalized L{\'e}vy measure $\alpha s^{-1-\alpha}ds/\Gamma(1-\alpha)$ tends vaguely to $0$ on $(0,\infty)$, and the scale-invariant candidate $s^{-1}ds$ for an ``index $0$'' L{\'e}vy measure is not one. The appropriate substitutes for a ``stable subordinator of index $0$'' are the extremal processes and the Dickman subordinator. Related renormalized limits, including the Kasahara extremal
process and the Dickman subordinator, are studied in the companion paper \cite{LiDickman}.
\end{rem}

The weak convergence of Proposition \ref{H-convergence} can be supplemented by a formal first-order expansion. Since $\xi/(1+\xi)^2$ is the Laplace transform of $(1-x)e^{-x}$, the Mittag--Leffler expansion (\ref{ML-rate}) formally corresponds to
$$
p_\alpha(x,t)=e^{-x}\bigl[1-\alpha(\log t+\gamma)(1-x)\bigr]+o(\alpha).
$$
The first-order coefficient has zero integral and changes sign at $x=1$. It therefore describes a signed redistribution of mass across the mean of the limiting exponential law, rather than a new probability density. The correction vanishes at $t=e^{-\gamma}$, matching the operator-level expansion in Corollary \ref{cor-rate-expansion}. Through the subordination formula (\ref{subordination-general}), this formal calculation is consistent with that operator-valued expansion. We do not pursue uniform density estimates here.

\subsection{The Laplacian}

In this subsection we specialize the abstract results of Section \ref{sec-abstract} to $A=c^2\Delta$ on the whole space $\mathbb{R}^n$ and on bounded domains, where the limit can be computed explicitly.

\begin{prop}\label{thm-Rn}
Let $c>0$ and let $A=c^2\Delta$ with its standard realization in $X=L^q(\mathbb{R}^n)$ ($1<q<\infty$) or $X=C_0(\mathbb{R}^n)$; in all these spaces $A$ generates a bounded analytic $C_0$-semigroup. Let $\phi\in X$ and let $f$ be as in Theorem \ref{thm-nonhomogeneous}. Let $u_\alpha$ be the solution of
\begin{equation*}
\begin{cases}
D_t^\alpha u(x,t)=c^2\Delta_x u(x,t)+f(x,t),  \quad  x\in\mathbb{R}^n,\ t>0,\\
u(x,0)=\phi(x).
\end{cases}
\end{equation*}
Then for every $t>0$,
\begin{equation*}
\lim_{\alpha\rightarrow0^+}u_\alpha(\cdot,t)=(1-c^2\Delta)^{-1}(\phi+f(\cdot,t)) \quad\text{in } X,
\end{equation*}
i.e. the limit $u(\cdot,t)$ is the unique solution of the elliptic equation
\begin{equation}\label{limit-elliptic}
u(x,t)-c^2\Delta_x u(x,t)=\phi(x)+f(x,t),\qquad x\in\mathbb{R}^n,
\end{equation}
in which $t>0$ merely plays the role of a parameter. Explicitly, the limit is given by convolution with the Bessel potential type kernel
\begin{equation}\label{Bessel}
(1-c^2\Delta)^{-1}\psi(x)=\int_{\mathbb{R}^n} G_c(x-y)\psi(y)\,dy,\qquad G_c(x)=\int_0^\infty e^{-s}(4\pi c^2 s)^{-n/2}e^{-\frac{|x|^2}{4c^2s}}\,ds.
\end{equation}
In dimension $n=1$ the kernel is elementary: $G_c(x)=e^{-|x|/c}/(2c)$.
\end{prop}
\begin{proof}
Theorems \ref{thm-homogeneous} and \ref{thm-nonhomogeneous} apply directly. The kernel representation (\ref{Bessel}) follows from $(1-c^2\Delta)^{-1}=\int_0^\infty e^{-s}e^{c^2s\Delta}\,ds$ and the explicit heat kernel. For $n=1$, the classical integral $\int_0^\infty s^{-1/2}e^{-s-a/s}ds=\sqrt{\pi}\,e^{-2\sqrt{a}}$ with $a=x^2/(4c^2)$ gives $G_c(x)=e^{-|x|/c}/(2c)$, and one checks directly that $(1-c^2\partial_x^2)G_c=\delta_0$.
\end{proof}

\begin{prop}
Let $\Omega\subset\mathbb{R}^n$ be a bounded domain with $C^{1,1}$ boundary, let $c>0$, and let $A=c^2\Delta_D$ be the Dirichlet Laplacian on $X=L^q(\Omega)$ ($1<q<\infty$), which generates a bounded analytic $C_0$-semigroup. Let $\phi$ and $f$ be as in Proposition \ref{thm-Rn} and let $u_\alpha$ be the solution of
\begin{equation*}
\begin{cases}
D_t^\alpha u(x,t)=c^2\Delta_x u(x,t)+f(x,t),  \quad  x\in\Omega,\ t>0,\\
u(x,t)=0,\quad x\in\partial\Omega,\ t>0,\\
u(x,0)=\phi(x).
\end{cases}
\end{equation*}
Then for every $t>0$, $u_\alpha(\cdot,t)$ converges in $L^q(\Omega)$ to the unique solution of the elliptic boundary value problem
\begin{equation}\label{limit-elliptic-bdd}
\begin{cases}
u(x,t)-c^2\Delta_x u(x,t)=\phi(x)+f(x,t),\quad x\in\Omega,\\
u(x,t)=0,\quad x\in\partial\Omega.
\end{cases}
\end{equation}
\end{prop}
\begin{proof}
This is a direct application of Theorems \ref{thm-homogeneous} and \ref{thm-nonhomogeneous} to $A=c^2\Delta_D$.
\end{proof}

\begin{rem}
(Degenerate elliptic equations.) The limit equations (\ref{limit-elliptic}) and (\ref{limit-elliptic-bdd}) are elliptic in the space variable and degenerate in time, the time variable entering only as a parameter; their well-posedness and regularity are classical. What our results add is the variational-in-$\alpha$ interpretation: these degenerate elliptic equations are exactly the $\alpha\rightarrow0^+$ limits of time fractional diffusion equations, in the same way as parabolic equations are the $\alpha\rightarrow1^-$ limits \cite{CT2023}.
\end{rem}

\subsection{Space-time fractional equations: comparison with the spatial limit \texorpdfstring{$\beta\rightarrow0^+$}{beta to 0+}}\label{sec-spacetime}

Expansions in the order also have a spatial counterpart. The logarithmic Laplacian arises as the first derivative of $(-\Delta)^s$ at $s=0$, and its higher-order versions provide finite Taylor expansions of the fractional Laplacian in its order. The expansions developed here concern the corresponding Caputo and Mittag-Leffler problems \cite{ChenWeth2019,ChenMOrder2024}.

For the spatial fractional Laplacian, the pointwise limit
\begin{equation*}
\lim_{\beta\rightarrow 0^+} (-\Delta)^\beta v(x)=v(x),
\end{equation*}
valid for instance for every $v$ in the Schwartz class by Fourier domination (constant functions being a counterexample to any unrestricted claim), has been widely studied (see e.g. \cite{Stinga2019}, and \cite[Remark 2.1, Proposition 2.1]{DaoudLaamri2022}): the spatial operator converges to the identity, in contrast with the Caputo time derivative, which converges to $\phi\mapsto\phi-\phi(0)$ (Proposition \ref{Caputo-limit-thm}(1)). The difference is explained on the Fourier/Laplace side: the symbol $|\xi|^{2\beta}$ of $(-\Delta)^\beta$ tends for every $\xi\neq0$ to $1$, the symbol of $I$ (at $\xi=0$ it equals $0$), while the Laplace symbol of $D_t^\alpha$ carries the initial-value term $\lambda^{\alpha-1}\phi(0)$, which survives as $\lambda^{-1}\phi(0)$.

We now study the joint limit for the space-time fractional equation
\begin{equation}\label{space-time-FDE}
\begin{cases}
D_t^\alpha u_{\alpha,\beta}(x,t)=-(-\Delta)^\beta u_{\alpha,\beta}(x,t),  \quad x\in\mathbb{R}^n,\ t>0, \\
u_{\alpha,\beta}(x,0)=\phi(x).
\end{cases}
\end{equation}

\begin{prop}\label{thm-spacetime}
Let $\phi\in L^2(\mathbb{R}^n)$ and let $u_{\alpha,\beta}$ be the solution of (\ref{space-time-FDE}) with $0<\alpha,\beta<1$.

(1) For fixed $\beta$, $\lim_{\alpha\rightarrow0^+}u_{\alpha,\beta}(\cdot,t)=(1+(-\Delta)^\beta)^{-1}\phi$ in $L^2(\mathbb{R}^n)$ for every $t>0$.

(2) As $(\alpha,\beta)\rightarrow(0,0)$ jointly, for every $t>0$,
\begin{equation}\label{joint-limit}
u_{\alpha,\beta}(\cdot,t)\longrightarrow \frac12\,\phi \qquad\text{in } L^2(\mathbb{R}^n).
\end{equation}
The iterated limits $\alpha\rightarrow0^+$ then $\beta\rightarrow0^+$, and $\beta\rightarrow0^+$ then $\alpha\rightarrow0^+$, give the same value $\phi/2$.
\end{prop}
\begin{proof}
Taking the Fourier transform in $x$, $\widetilde{u}_{\alpha,\beta}(\xi,t)=E_\alpha(-|\xi|^{2\beta}t^\alpha)\widetilde\phi(\xi)$.

(1) $-(-\Delta)^\beta$ generates the bounded symmetric $2\beta$-stable semigroup on $L^2$, so Theorem \ref{thm-homogeneous} applies.

(2) Write $s_{\alpha,\beta}:=|\xi|^{2\beta}t^\alpha$; for $\xi\neq0$. Then $s_{\alpha,\beta}\rightarrow1$ as $(\alpha,\beta)\rightarrow(0,0)$. By (\ref{ML-Laplace}), $E_\alpha(-s)=\int_0^\infty e^{-sx}p_\alpha(x,1)\,dx$, and $p_\alpha(\cdot,1)$ converges weakly to $e^{-x}dx$ by Proposition \ref{H-convergence}. Hence, writing
$$
\left|E_\alpha(-s_{\alpha,\beta})-\frac12\right|\leq \left|E_\alpha(-s_{\alpha,\beta})-E_\alpha(-1)\right|+\left|E_\alpha(-1)-\frac12\right|,
$$
the second term tends to $0$ by Lemma \ref{ML-limit}. For the first term, once $s_{\alpha,\beta}\geq1/2$, the inequality $|e^{-sx}-e^{-x}|\leq |s-1|xe^{-x/2}$ for $s\geq1/2$ and $x\geq0$ gives the bound
$$
|s_{\alpha,\beta}-1|\int_0^\infty x e^{-x/2}p_\alpha(x,1)\,dx\leq (2/e)|s_{\alpha,\beta}-1|\rightarrow0,
$$
since $xe^{-x/2}\leq2/e$. Thus $\widetilde u_{\alpha,\beta}(\xi,t)\rightarrow\widetilde\phi(\xi)/2$ for a.e. $\xi$, with $|\widetilde u_{\alpha,\beta}(\xi,t)|\leq|\widetilde\phi(\xi)|$ (as $|E_\alpha(-s)|\leq1$ for $s\geq0$ by complete monotonicity), and dominated convergence gives (\ref{joint-limit}).

The iterated limits are read off from the solution formula. If $\alpha\rightarrow0^+$ first, (1) gives $u_{\alpha,\beta}(\cdot,t)\rightarrow(1+(-\Delta)^\beta)^{-1}\phi$, and since $(1+|\xi|^{2\beta})^{-1}\widetilde\phi(\xi)\rightarrow\frac12\widetilde\phi(\xi)$ for a.e. $\xi$ as $\beta\rightarrow0^+$, with $|(1+|\xi|^{2\beta})^{-1}\widetilde\phi(\xi)|\leq|\widetilde\phi(\xi)|$, the second limit is $\phi/2$ by dominated convergence. If $\beta\rightarrow0^+$ first, then for every fixed $\alpha$, $\widetilde u_{\alpha,\beta}(\xi,t)\rightarrow E_\alpha(-t^\alpha)\widetilde\phi(\xi)$ for a.e. $\xi$, bounded in modulus by $|\widetilde\phi(\xi)|$. Hence $u_{\alpha,\beta}(\cdot,t)\rightarrow E_\alpha(-t^\alpha)\phi$ in $L^2$. This limit is the solution of $D_t^\alpha u=-u$, $u(0)=\phi$, in accordance with the spatial limit $-(-\Delta)^\beta\rightarrow-I$. The subsequent limit $\alpha\rightarrow0^+$ is the case $\beta=0$ of the computation above, where $s_{\alpha,0}=t^\alpha\rightarrow1$, and equals $\phi/2$.
\end{proof}

\begin{rem}\label{rem-spacetime-stochastic}
Stochastically, $u_{\alpha,\beta}(x,t)=\mathbb{E}[\phi(x+S_{2\beta}(H_\alpha(t)))]$, where $S_{2\beta}$ is a symmetric $2\beta$-stable process independent of $H_\alpha$. For a probability density $\phi\geq0$ belonging to $L^2$, the joint limit $\phi/2$ has the following mass interpretation. As $\beta\rightarrow0^+$, the characteristic function $\mathbb{E}e^{i\xi S_{2\beta}(s)}=e^{-s|\xi|^{2\beta}}\rightarrow e^{-s}$ ($\xi\neq0$) shows that the process $S_{2\beta}(s)$ retains only the mass $e^{-s}$ at the origin, while the remaining mass $1-e^{-s}$ escapes to infinity. The retained mass is then obtained by averaging $e^{-s}$ over the limiting exponential clock, as explained below. Since each $u_{\alpha,\beta}(\cdot,t)$ has total mass one, the convergence in (\ref{joint-limit}) is in $L^2$, not in $L^1$. This degeneration parallels the escape of mass of the stable subordinator noted in Remark \ref{rem-p-alpha}(3).

Both degenerations are in fact the same phenomenon, seen once through the inverse and once directly. Write $S_{2\beta}(s)\overset{d}{=}B(2Y_\beta(s))$, where $B$ is a Brownian motion and $Y_\beta$ a $\beta$-stable subordinator, $\mathbb{E}e^{-\lambda Y_\beta(s)}=e^{-s\lambda^\beta}$, independent of $B$ and of $H_\alpha$; then $u_{\alpha,\beta}(x,t)=\mathbb{E}[\phi(x+B(2Y_\beta(H_\alpha(t))))]$, a Brownian motion under two nested stable-type clocks. As $\beta\rightarrow0^+$ the Laplace exponent $\lambda^\beta$ tends to $1$, so $Y_\beta(s)\Rightarrow e^{-s}\delta_0+(1-e^{-s})\delta_\infty$ on $[0,\infty]$: $Y_\beta$ converges to the degenerate subordinator $Y_0$ with exponent $\psi\equiv1$, which stays at $0$ up to an $\mathrm{Exp}(1)$ killing time $\zeta$ and then jumps to $+\infty$. On the time side, the limit $\tau$ of the inverse clocks $H_\alpha(t)$ (Proposition \ref{H-convergence}) is the inverse of an independent copy of the same degenerate subordinator, that is, its killing time, hence an $\mathrm{Exp}(1)$ random variable independent of $\zeta$. The limiting object $x+B(2Y_0(\tau))$ therefore equals $x$ on $\{\tau<\zeta\}$ and has escaped to infinity on $\{\tau\geq\zeta\}$, and the retained mass is
\begin{equation*}
\mathbb{P}(\tau<\zeta)=\int_0^\infty e^{-t}\,e^{-t}\,dt=\frac12,
\end{equation*}
the competition of two independent exponential clocks. This explains at once the factor $1/2$, the $t$-independence of the limit (the limiting clock $\tau$ does not depend on $t$), and the symmetry of the two iterated limits (interchanging the order interchanges the roles of $\tau$ and $\zeta$). On the Fourier--Laplace side the same two unit killing rates appear in the transform $\lambda^{\alpha-1}\widetilde\phi(\xi)/(\lambda^\alpha+|\xi|^{2\beta})\rightarrow\widetilde\phi(\xi)/(2\lambda)$ of the solution.

\end{rem}

\subsection{The difference operator: the fractional Poisson process}

A basic instance of (\ref{FCP}) with a bounded operator is the governing equation of the fractional Poisson process. The state probabilities $p_k(t)=\mathbb{P}(N_\alpha(t)=k)$ of the fractional Poisson process $N_\alpha(t)$ of rate $\lambda>0$ are determined by the fractional difference-differential equation
\begin{equation}\label{FPP-eq}
D_t^\alpha p_k(t)=-\lambda\bigl(p_k(t)-p_{k-1}(t)\bigr),\quad t>0,\ k\geq0,\qquad p_k(0)=\delta_{k0},
\end{equation}
with $p_{-1}(t)\equiv0$; see \cite{BeghinOrsingher2009FPP, Orsingher2013FPP}, and also \cite{LizamaRebolledo2018} for a semigroup approach. This is precisely (\ref{FCP}) with $X=\ell^1(\mathbb{N}_0)$, $\phi=e_0=(\delta_{k0})_{k\geq0}$ and the difference operator
\begin{equation}\label{difference-operator}
A=-\lambda(I-R),\qquad (Rp)_k=p_{k-1}\mathbf{1}_{\{k\geq1\}},
\end{equation}
a bounded operator on $X$ generating the uniformly continuous forward semigroup $T(s)=e^{-\lambda s(I-R)}$ of the standard Poisson process of rate $\lambda$ (a positive, mass-preserving contraction semigroup on $\ell^1$). The subordination formula (\ref{subordination2}) therefore gives
\begin{equation}\label{FPP-subordination}
p_k(t)=\mathbb{E}\left[e^{-\lambda H_\alpha(t)}\frac{(\lambda H_\alpha(t))^k}{k!}\right],
\end{equation}
that is, $N_\alpha(t)\stackrel{d}{=}N_1(H_\alpha(t))$ with $N_1$ a standard Poisson process of rate $\lambda$ independent of $H_\alpha$; the expectation in (\ref{FPP-subordination}) equals the classical formula $\mathbb{P}(N_\alpha(t)=k)=\frac{(\lambda t^\alpha)^k}{k!}E_\alpha^{(k)}(-\lambda t^\alpha)$, see \cite{MNV2011}. The abstract representation (\ref{subordination2}) thus identifies the two definitions of the fractional Poisson process used in the literature: the equation definition (\ref{FPP-eq}) and the renewal/time-change definition with Mittag-Leffler interarrival times.

\begin{prop}
For every $t>0$, $N_\alpha(t)$ converges in total variation as $\alpha\rightarrow0^+$ to the geometric random variable
\begin{equation*}
\lim_{\alpha\rightarrow0^+}\mathbb{P}(N_\alpha(t)=k)=\frac{\lambda^k}{(1+\lambda)^{k+1}},\qquad k=0,1,2,\dots
\end{equation*}
that is, to $N_1(\tau)$ with $\tau\sim\mathrm{Exp}(1)$ independent of $N_1$. The limit is again independent of $t$.
\end{prop}
\begin{proof}
%We give two proofs. First,
Theorem \ref{thm-homogeneous} applied to the bounded operator (\ref{difference-operator}) gives $p(t)\rightarrow q:=(1-A)^{-1}e_0$ in $\ell^1$, i.e.\ in total variation, and the equation $(1-A)q=e_0$ reads $(1+\lambda)q_k-\lambda q_{k-1}=\delta_{k0}$, hence $q_k=\lambda^k/(1+\lambda)^{k+1}$.
%Second, by (\ref{FPP-subordination}) and Proposition \ref{H-convergence}, the bounded continuous function $s\mapsto e^{-\lambda s}(\lambda s)^k/k!$ gives
%$$
%\mathbb{P}(N_\alpha(t)=k)=\mathbb{E}\left[e^{-\lambda H_\alpha(t)}\frac{(\lambda H_\alpha(t))^k}{k!}\right]\longrightarrow
%\int_0^\infty e^{-\lambda s}\frac{(\lambda s)^k}{k!}e^{-s}\,ds=\frac{\lambda^k}{(1+\lambda)^{k+1}}.
%$$
%Since the limiting probabilities sum to one, Scheff\'e's lemma upgrades the second argument to convergence in total variation as well.
\end{proof}

\begin{rem}
(1) The geometric distribution is the discrete counterpart of the exponential limit of $H_\alpha(t)$: the fractional Poisson process degenerates to a mixed Poisson random variable $N_1(\tau)$. Consistently, the survival function of the first interarrival time satisfies $\mathbb{P}(T_1^\alpha>t)=E_\alpha(-\lambda t^\alpha)\rightarrow(1+\lambda)^{-1}$, which is $\mathbb{P}(N_1(\tau)=0)$.

(2) The same computation applies to counting processes time-changed by $H_\alpha$ generally: any event probability that is a bounded continuous function of the clock converges to its exponential-clock counterpart.
\end{rem}

The space-fractional Poisson process of \cite{OrsingherPolito2012} fractionalizes instead the difference operator, replacing $I-R$ in (\ref{FPP-eq}) by its fractional power $(I-R)^\beta=\sum_{j\geq0}(-1)^j\binom{\beta}{j}R^j$ (the fractional difference operator of time series analysis): the governing equation becomes $D_t^\alpha p=-\lambda^\beta(I-R)^\beta p=-(-A)^\beta p$, the sequence-space analogue of the space-time equation (\ref{space-time-FDE}); its $\beta\rightarrow0^+$ behavior parallels Section \ref{sec-spacetime} exactly (Proposition \ref{prop-spacetime-FPP} below).

\begin{prop}\label{prop-spacetime-FPP}
Let $N_{\alpha,\beta}(t)$ be the space-time fractional Poisson process of \cite{OrsingherPolito2012}, that is, the counting process with probability generating function
\begin{equation}\label{pgf-spacetime-FPP}
\mathbb{E}\left[u^{N_{\alpha,\beta}(t)}\right]=E_\alpha\left(-\lambda^\beta(1-u)^\beta t^\alpha\right),\qquad |u|\leq1,
\end{equation}
with the principal branch of $(1-u)^\beta$ and the convention $0^\beta:=0$.
Then, as $(\alpha,\beta)\rightarrow(0,0)$, for every $t>0$,
\begin{equation*}
\mathbb{P}(N_{\alpha,\beta}(t)=0)\longrightarrow\frac12,\qquad \mathbb{P}(N_{\alpha,\beta}(t)\geq k)\longrightarrow\frac12\quad\text{for every }k\geq1;
\end{equation*}
equivalently, on the compactified state space $\mathbb{N}_0\cup\{\infty\}$, $N_{\alpha,\beta}(t)\Rightarrow\frac12\delta_0+\frac12\delta_\infty$: half of the mass remains at $0$ and half escapes to infinity, exactly as in Proposition \ref{thm-spacetime} and Remark \ref{rem-spacetime-stochastic}.
\end{prop}
\begin{proof}
For fixed $u\in[0,1)$ one has $\lambda^\beta(1-u)^\beta t^\alpha\rightarrow1$, and the argument of Proposition \ref{thm-spacetime} (2) (writing $E_\alpha(-s)=\int_0^\infty e^{-sx}p_\alpha(x,1)\,dx$ and using $xe^{-x/2}\leq2/e$) gives $E_\alpha(-\lambda^\beta(1-u)^\beta t^\alpha)\rightarrow\frac12$. In particular $\mathbb{P}(N_{\alpha,\beta}(t)=0)=E_\alpha(-\lambda^\beta t^\alpha)\rightarrow\frac12$, the case $u=0$. For $k\geq1$ and any fixed $u\in(0,1)$,
$$
\mathbb{P}(1\leq N_{\alpha,\beta}(t)\leq k-1)\leq u^{-(k-1)}\left(\mathbb{E}\left[u^{N_{\alpha,\beta}(t)}\right]-\mathbb{P}(N_{\alpha,\beta}(t)=0)\right)\longrightarrow u^{-(k-1)}\left(\frac12-\frac12\right)=0,
$$
hence $\mathbb{P}(N_{\alpha,\beta}(t)\geq k)=1-\mathbb{P}(N_{\alpha,\beta}(t)=0)-\mathbb{P}(1\leq N_{\alpha,\beta}(t)\leq k-1)\rightarrow\frac12$. The distributional assertion follows, the limit law on $\mathbb{N}_0\cup\{\infty\}$ being determined by its point masses. Stochastically, (\ref{pgf-spacetime-FPP}) is the probability generating function of $N_1(S_\beta(H_\alpha(t)))$, with $S_\beta$ a $\beta$-stable subordinator independent of $H_\alpha$. Its mass-escape mechanism is exactly the two-clock mechanism of Remark \ref{rem-spacetime-stochastic}.
\end{proof}

\subsection{Potential theory and the stochastic interpretation of the limit}\label{sec-potential}

\subsubsection{The limit solution as a $\lambda$-potential}

Let $X(t)$ be a Markov process with semigroup $(T(t))_{t\geq0}$ and generator $A$ (for instance Brownian motion, for which $A=\frac12\Delta$). In probabilistic potential theory, the resolvent operators
\begin{equation}\label{lambda-potential}
(\lambda-A)^{-1}f=\int_0^\infty e^{-\lambda s}T(s)f\,ds=\frac1\lambda\,\mathbb{E}\left[T(\tau_\lambda)f\right],\qquad \lambda>0,
\end{equation}
where $\tau_\lambda$ is an exponential random variable with parameter $\lambda$ independent of $X$, are called the $\lambda$-potentials of $X(t)$. Probabilistically, $\lambda(\lambda-A)^{-1}f$ is the expected value of $f(X(s))$ sampled at the independent exponential time $\tau_\lambda$, see e.g. \cite{Schilling2012}.

The results of Section \ref{sec-abstract} can now be restated in these terms, with an additional scale parameter.

\begin{prop}
Suppose $A$ generates a bounded $C_0$-semigroup $(T(t))_{t\geq0}$ on $X$. Let $\lambda>0$ and let $u_\alpha$ be the mild solution of the rescaled fractional Cauchy problem
\begin{equation}\label{lambda-FCP}
\begin{cases}
\lambda D_t^\alpha u_\alpha(t)=Au_\alpha(t), \quad t>0, \\
\lambda u_\alpha(0)=f.
\end{cases}
\end{equation}
Then for every $t>0$,
\begin{equation}\label{lambda-limit}
\lim_{\alpha\rightarrow0^+}u_\alpha(t)=(\lambda-A)^{-1}f,
\end{equation}
uniformly on compact subsets of $(0,\infty)$. In particular, for $\lambda=1$ one recovers Theorem \ref{thm-homogeneous}.
\end{prop}
\begin{proof}
We give two proofs. First, (\ref{lambda-FCP}) is equivalent to $D_t^\alpha u_\alpha=(A/\lambda)u_\alpha$, $u_\alpha(0)=f/\lambda$, and $A/\lambda$ generates the bounded semigroup $T(s/\lambda)$; Theorem \ref{thm-homogeneous} yields
$$
u_\alpha(t)\longrightarrow\left(1-\frac{A}{\lambda}\right)^{-1}\frac{f}{\lambda}=(\lambda-A)^{-1}f.
$$
Second, the function $s\mapsto T(s/\lambda)f/\lambda$ is bounded and continuous. Applying Lemma \ref{lem-measure-weak} to this function in the subordination formula, and then substituting $s=\lambda r$, gives
$$
u_\alpha(t)=\int_0^\infty p_\alpha(s,t)\,T\!\left(\frac{s}{\lambda}\right)\frac{f}{\lambda}\,ds
\longrightarrow\int_0^\infty e^{-s}T\!\left(\frac{s}{\lambda}\right)\frac{f}{\lambda}\,ds
=\int_0^\infty e^{-\lambda r}T(r)f\,dr=(\lambda-A)^{-1}f.
$$
\end{proof}

\begin{cor}
Let $X(t)$ be a Feller process on a locally compact separable metric space $E$, with semigroup $T(s)$ on $C_0(E)$ generated by $A$, and independent of the inverse stable subordinator $H_\alpha$. Then for every bounded continuous $f$, every $x$ and every $t,\lambda>0$,
\begin{equation*}
\lim_{\alpha\rightarrow0^+}\mathbb{E}_x\left[f\bigl(X(H_\alpha(t)/\lambda)\bigr)\right]=\mathbb{E}_x\left[f\bigl(X(\tau_\lambda)\bigr)\right],
\end{equation*}
where $\tau_\lambda$ is exponential of rate $\lambda$, independent of $X$; for $f\in C_0(E)$ the limit equals $\lambda(\lambda-A)^{-1}f(x)$.
\end{cor}
\begin{proof}
Conditioning on the independent clock gives
$$
\mathbb E_x f\bigl(X(H_\alpha(t)/\lambda)\bigr)
=\int_0^\infty F(s)p_\alpha(s,t)\,ds,
\qquad
F(s):=\mathbb E_x f\bigl(X(s/\lambda)\bigr).
$$
The function $F$ is bounded and continuous by the stochastic continuity of the Feller process. Hence Proposition \ref{H-convergence} gives
$$
\int_0^\infty F(s)p_\alpha(s,t)\,ds
\longrightarrow\int_0^\infty e^{-s}F(s)\,ds.
$$
The limiting integral is $\mathbb E_x f(X(\tau_\lambda))$ by the change of variables $s=\lambda r$. For $f\in C_0(E)$ it equals $\lambda(\lambda-A)^{-1}f(x)$ by (\ref{lambda-potential}).
\end{proof}

\begin{rem}
For each fixed $\alpha>0$, the process $X(H_\alpha(t))$ is generally non-Markov because its clock $H_\alpha(t)$ has memory. As $\alpha\rightarrow0^+$, this memory collapses and only the exponential clock survives. The convergence asserted here concerns one-time marginals, not path-level convergence. Moreover, (\ref{lambda-FCP})--(\ref{lambda-limit}) provide a fractional approximation of resolvents: for small $\alpha$, the value $u_\alpha(t)$ at any $t>0$ approximates $(\lambda-A)^{-1}f$.
\end{rem}

\subsubsection{Relation with the Feynman--Kac formula}

The classical Feynman--Kac formula represents elliptic resolvents by Brownian functionals \cite[pp. 110--112]{Evans2013}. In the present whole-space setting, let $c>0$, let $\varphi$ and $f(\cdot,t)$ be bounded continuous, and let $u$ be the bounded distributional solution of
\begin{equation}\label{limit-eq-FK}
c\,\bigl(u(x,t)-\varphi(x)\bigr)=\frac12\Delta_x u(x,t)+f(x,t),\qquad x\in\mathbb{R}^n,\ t>0.
\end{equation}
Then
\begin{equation*}
u(x,t)=\mathbb{E}\left[\varphi\bigl(x+B(\tau_c)\bigr)\right]+\mathbb{E}\left[\int_0^{\tau_c} f\bigl(x+B(r),t\bigr)\,dr\right],
\end{equation*}
where $B$ is a standard Brownian motion and $\tau_c$ an independent $\mathrm{Exp}(c)$ random variable. This is the probabilistic form of
$$
u(\cdot,t)=\left(c-\frac12\Delta\right)^{-1}(c\varphi+f(\cdot,t))
=c\int_0^\infty e^{-cr}T(r)\varphi\,dr+\int_0^\infty e^{-cr}T(r)f(\cdot,t)\,dr,
$$
with $T(r)$ the Brownian semigroup. Whenever the rescaled problem satisfies the hypotheses of Theorem \ref{thm-nonhomogeneous}, this $u$ is the $\alpha\rightarrow0^+$ limit of $cD_t^\alpha u_\alpha=\frac12\Delta_x u_\alpha+f$ with $u_\alpha(\cdot,0)=\varphi$.

For comparison, regarding $t$ as an additional space variable gives the different elliptic problem
\begin{equation}\label{halfspace-problem}
-\frac12\Delta_{x,t}u+cu=f\quad\text{in }\mathbb R^n\times(0,\infty),
\qquad u(x,0)=\varphi(x),
\end{equation}
in which $\Delta_{x,t}$ is the full Laplacian. To relate the two problems, set
$$
\mathcal L:=\sqrt{2c-\Delta_x}.
$$
With $f(\cdot,t)=g$, they read respectively
$$
\mathcal L^2u=2c\varphi+2g,
\qquad
(\mathcal L^2-\partial_t^2)u=2g\quad\text{in }\mathbb R^n\times(0,\infty),\quad u(\cdot,0)=\varphi.
$$
Thus the extension applies $e^{-t\mathcal L}$, while the limit equation applies $2c\mathcal L^{-2}$.

Fix $c>0$ and $\varphi,g\in C_b(\mathbb{R}^n)$. Let
$$
u_A(x):=\Bigl(c-\tfrac12\Delta_x\Bigr)^{-1}(c\varphi+g)(x)
$$
be the solution of (\ref{limit-eq-FK}) with $f(\cdot,t)=g$.
Let
$$
u_B(x,t)=\mathbb{E}\Big[\varphi\bigl(x+B_1(\tau_t)\bigr)e^{-c\tau_t}\Big]+\mathbb{E}\left(\int_0^{\tau_t}g\bigl(x+B_1(r)\bigr)e^{-cr}\,dr\right)
$$
be the bounded distributional solution of (\ref{halfspace-problem}) with pointwise boundary values, where $\tau_t$ is the first hitting time of $0$ by $t+B_2$. We denote the $\varphi$-part of $u_B$ by $u_B^{(0)}$. Finally, for $\varepsilon>0$ let $u_\varepsilon$ be the bounded distributional solution, with pointwise boundary values, of
\begin{equation}\label{anisotropic}
-\frac12\Delta_x u_\varepsilon-\frac{\varepsilon}{2}\,\partial_t^2 u_\varepsilon+cu_\varepsilon=g\quad\text{in }\mathbb{R}^n\times(0,\infty),\qquad u_\varepsilon(\cdot,0)=\varphi.
\end{equation}

\begin{thm}
With the above notation, the following statements hold:

{\rm (i)} (Anisotropic degeneration.) For every fixed $t>0$, $u_\varepsilon(x,t)\rightarrow(c-\frac12\Delta_x)^{-1}g(x)$ pointwise as $\varepsilon\rightarrow0^+$. More precisely, $u_\varepsilon=u_\varepsilon^{(0)}+u_\varepsilon^{(g)}$, where
$$
u_\varepsilon^{(0)}(x,t)=e^{-t\sqrt{(2c-\Delta_x)/\varepsilon}}\varphi(x)=\mathbb{E}\Big[\varphi\bigl(x+B_1\bigl(\tau_t^{(\varepsilon)}\bigr)\bigr)e^{-c\tau_t^{(\varepsilon)}}\Big],
$$
$$
\sup_x|u_\varepsilon^{(0)}(x,t)|\leq\|\varphi\|_\infty\,e^{-t\sqrt{2c/\varepsilon}},
$$
with $\tau_t^{(\varepsilon)}$ the hitting time of $0$ by $t+\sqrt\varepsilon\,B_2$ and $\mu_t^{(\varepsilon)}$ its law, and with $u_\varepsilon^{(g)}(x,t)\rightarrow(c-\frac12\Delta_x)^{-1}g(x)$. Thus the boundary datum survives only in a boundary layer of width $\sqrt\varepsilon$.

{\rm (ii)} (Moment identity.) The laws $\mu_t$ of $\tau_t$ satisfy
\begin{equation}\label{exit-moment}
\int_0^\infty 2t\,\mu_t(ds)\,dt=ds\qquad\text{on }(0,\infty).
\end{equation}
Equivalently,
\begin{equation}\label{exit-moment-discounted}
ce^{-cs}\,ds=\int_0^\infty 2ct\,e^{-cs}\,\mu_t(ds)\,dt\qquad\text{on }(0,\infty),
\end{equation}
so the exponential clock of the $\alpha\rightarrow0^+$ limit is the $2ct\,dt$-weighted average of the discounted half-space exit-time clocks. The $\varphi$-parts of the two solutions are linked by
\begin{equation}\label{halfspace-moment}
u_A^{(\varphi)}(x):=\mathbb{E}\big[\varphi\bigl(x+B(\tau_c)\bigr)\big]=2c\int_0^\infty t\,u_B^{(0)}(x,t)\,dt,
\end{equation}
that is, $2c\,\mathcal L^{-2}=2c\int_0^\infty t\,e^{-t\mathcal L}dt$.
\end{thm}
\begin{proof}
The hitting-time Laplace transform is
$$
\mathbb E[e^{-\lambda\tau_t}]=e^{-t\sqrt{2\lambda}},
\qquad
\mathbb E[e^{-\lambda\tau_t^{(\varepsilon)}}]=e^{-t\sqrt{2\lambda/\varepsilon}}.
$$
In particular, $\tau_t^{(\varepsilon)}\rightarrow\infty$ in probability as $\varepsilon\rightarrow0^+$.

(i) The Feynman--Kac formula applied to (\ref{anisotropic}) gives the stated decomposition. The estimate for $u_\varepsilon^{(0)}$ follows from
$|u_\varepsilon^{(0)}(x,t)|\leq\|\varphi\|_\infty\mathbb E[e^{-c\tau_t^{(\varepsilon)}}]$. For the forcing part, dominated convergence in
$$
\mathbb E\left(\int_0^{\tau_t^{(\varepsilon)}}g(x+B_1(r))e^{-cr}\,dr\right)
=\int_0^\infty e^{-cr}T(r)g(x)\,\mathbb P(\tau_t^{(\varepsilon)}>r)\,dr
$$
gives $u_\varepsilon^{(g)}(x,t)\rightarrow\int_0^\infty e^{-cr}T(r)g(x)\,dr=(c-\frac12\Delta_x)^{-1}g(x)$. The operator form of $u_\varepsilon^{(0)}$ follows by taking Fourier symbols in the subordination representation
$u_\varepsilon^{(0)}=\int_0^\infty T(s)\varphi\,e^{-cs}\,\mu_t^{(\varepsilon)}(ds)$.

(ii) Taking Laplace transforms in $s$ gives
$$
\int_0^\infty e^{-\lambda s}\int_0^\infty2t\,\mu_t(ds)\,dt
=\int_0^\infty2t\,e^{-t\sqrt{2\lambda}}\,dt
=\frac1\lambda,
$$
which proves (\ref{exit-moment}) and hence (\ref{exit-moment-discounted}). Subordination gives
$$
u_B^{(0)}(x,t)=\int_0^\infty T(s)\varphi(x)\,e^{-cs}\,\mu_t(ds).
$$
Multiplying by $2ct$ and integrating in $t$, using Fubini's theorem and (\ref{exit-moment}), yields
$$
2c\int_0^\infty t\,u_B^{(0)}(x,t)\,dt
=\int_0^\infty T(s)\varphi(x)\,ce^{-cs}\,ds
=\mathbb E[\varphi(x+B(\tau_c))].
$$
The final operator identity is $\mathcal L^{-2}=\int_0^\infty te^{-t\mathcal L}\,dt$.
\end{proof}

\begin{rem}
The two degenerations handle the datum differently: as $\alpha\rightarrow0^+$, the Caputo initial datum becomes the source term $c\varphi$, whereas as $\varepsilon\rightarrow0^+$ the elliptic boundary datum is confined to a layer of width $\sqrt\varepsilon$. The latter mechanism belongs to the classical theory of anisotropic singular perturbations \cite{Chipot2007,ChipotGuesmia2009,Lions1973,Ogabi2016}. The identity (\ref{halfspace-moment}) gives the explicit relation between the two constructions in the present half-space setting.
\end{rem}

\noindent{\bf Data availability.}
No datasets were generated or analyzed during this study.

\section*{Declaration of generative AI and AI-assisted technologies in the writing process}
During the preparation of this work, the author used Kimi K3 to assist with mathematical checking and English-language editing. After using this tool, the author reviewed and edited the content as needed and takes full responsibility for the content of the publication.


\begin{thebibliography}{99}

\bibitem{Apelblat2020} A. Apelblat, Differentiation of the Mittag-Leffler functions with respect to parameters in the Laplace transform approach, Mathematics {\bf 8} (2020), no. 5, 657.

\bibitem{Applebaum2009} D. Applebaum, L{\'e}vy Processes and Stochastic Calculus, 2nd ed., Cambridge Studies in Advanced Mathematics {\bf 116}, Cambridge University Press, Cambridge, 2009.

\bibitem{ABHN2011} W. Arendt, C. J. K. Batty, M. Hieber and F. Neubrander, Vector-valued Laplace transforms and Cauchy problems, 2nd ed., Monographs in Mathematics {\bf 96}, Birkh\"auser/Springer, Basel, 2011.

\bibitem{BM2001} B. Baeumer and M. M. Meerschaert, Stochastic solutions for fractional Cauchy problems, Fract. Calc. Appl. Anal. {\bf 4} (2001), no. 4, 481--500.

\bibitem{BMN2009} B. Baeumer, M. M. Meerschaert and E. Nane, Brownian subordinators and fractional Cauchy problems, Trans. Amer. Math. Soc. {\bf 361} (2009), no. 7, 3915--3930.

\bibitem{Baj} E. G. Bajlekova, Fractional evolution equations in Banach spaces, PhD thesis, Department of Mathematics, Eindhoven University of Technology, 2001.

\bibitem{BeghinOrsingher2009FPP} L. Beghin and E. Orsingher, Fractional Poisson processes and related planar random motions, Electron. J. Probab. {\bf 14} (2009), no. 61, 1790--1826.

\bibitem{Beghin2009} L. Beghin and E. Orsingher, Iterated elastic Brownian motions and fractional diffusion equations,
 Stochastic Process. Appl. {\bf 119} (2009), no. 6, 1975--2003.

\bibitem{CT2023} T. Caraballo and N. H. Tuan, New results for convergence problem of fractional diffusion equations when order approach to $1^-$, Differential Integral Equations {\bf 36} (2023), no. 5--6, 491--516.

% Temporarily retained for possible later discussion of iterated time changes.
% \bibitem{Casse2016} J. Casse and J.-F. Marckert,
%  Processes iterated ad libitum, Stochastic Process. Appl. {\bf 126} (2016), no. 11, 3353--3376.

\bibitem{ChenMOrder2024} H. Chen, On $m$-order logarithmic Laplacians and related properties, arXiv:2307.06198v4, 2024.

\bibitem{ChenWeth2019} H. Chen and T. Weth, The Dirichlet problem for the logarithmic Laplacian, Comm. Partial Differential Equations {\bf 44} (2019), no. 11, 1100--1139.

\bibitem{CNYY2009} J. Cheng, J. Nakagawa, M. Yamamoto and T. Yamazaki, Uniqueness in an inverse problem for a one-dimensional fractional diffusion equation, Inverse Problems {\bf 25} (2009), no. 11, 115002.

\bibitem{Chipot2007} M. Chipot, On some anisotropic singular perturbation problems, Asymptot. Anal. {\bf 55} (2007), no. 1--2, 125--144.

\bibitem{ChipotGuesmia2009} M. Chipot and S. Guesmia, On the asymptotic behavior of elliptic, anisotropic singular perturbations problems, Commun. Pure Appl. Anal. {\bf 8} (2009), no. 1, 179--193.

% \bibitem{Curien2014} N. Curien and T. Konstantopoulos,
%  Iterating Brownian motions, ad libitum. J. Theoret. Probab. {\bf 27} (2014), no. 2, 433--448.

\bibitem{DaoudLaamri2022} M. Daoud and E. H. Laamri, Fractional Laplacians: a short survey, Discrete Contin. Dyn. Syst. Ser. S {\bf 15} (2022), no. 1, 95--116.

\bibitem{Evans2013} L. C. Evans, An Introduction to Stochastic Differential Equations, American Mathematical Society, Providence, RI, 2013.

\bibitem{GKMR2020} R. Gorenflo, A. A. Kilbas, F. Mainardi and S. V. Rogosin, Mittag-Leffler Functions, Related Topics and Applications, 2nd ed., Springer Monographs in Mathematics, Springer, Berlin, 2020.

\bibitem{Hardy1949} G. H. Hardy, Divergent Series, Oxford, at the Clarendon Press, 1949.

\bibitem{HP1957} E. Hille and R. S. Phillips, Functional Analysis and Semi-Groups, revised ed., American Mathematical Society Colloquium Publications {\bf 31}, Amer. Math. Soc., Providence, RI, 1957.

\bibitem{LiDickman} C.-G. Li, Evolution equations driven by the Dickman derivative: continuity in the Laplace exponent, time-domain expansions, and the Kasahara scale, preprint, 2026.

\bibitem{LiChenLi2010} M. Li, C. Chen and F.-B. Li, On fractional powers of generators of fractional resolvent families, J. Funct. Anal. {\bf 259} (2010), no. 10, 2702--2726.

\bibitem{LiZheng2004} M. Li and Q. Zheng, On spectral inclusions and approximations of $\alpha$-times resolvent families, Semigroup Forum {\bf 69} (2004), no. 3, 356--368.

\bibitem{Lions1973} J. L. Lions, Perturbations singuli\`eres dans les probl\`emes aux limites et en contr\^ole optimal, Lecture Notes in Mathematics {\bf 323}, Springer-Verlag, Berlin, 1973.

\bibitem{LizamaRebolledo2018} C. Lizama and R. Rebolledo, A semigroup approach to fractional Poisson processes, Complex Anal. Oper. Theory {\bf 12} (2018), no. 3, 777--785.

\bibitem{MNV2011} M. M. Meerschaert, E. Nane and P. Vellaisamy, The fractional Poisson process and the inverse stable subordinator, Electron. J. Probab. {\bf 16} (2011), no. 59, 1600--1620.

\bibitem{Meerschaert2012book} M. M. Meerschaert and A. Sikorskii, Stochastic Models for Fractional Calculus,
De Gruyter Studies in Mathematics {\bf 43}, De Gruyter, Berlin, 2012.

\bibitem{Ogabi2016} C. Ogabi, On the $L^p$-theory of anisotropic singular perturbations of elliptic problems, Commun. Pure Appl. Anal. {\bf 15} (2016), no. 4, 1157--1178.

\bibitem{NIST} F. W. J. Olver, D. W. Lozier, R. F. Boisvert and C. W. Clark (eds.), NIST Handbook of Mathematical Functions, Cambridge University Press, Cambridge, 2010. Online companion: Digital Library of Mathematical Functions, \texttt{https://dlmf.nist.gov/}, Sections 5.7 and 5.9.

\bibitem{Orsingher2013FPP} E. Orsingher, Fractional Poisson processes, Sci. Math. Jpn. {\bf 76} (2013), no. 1, 139--145.

\bibitem{Orsingher2009} E. Orsingher and L. Beghin, Fractional diffusion equations and processes with randomly varying time,
Ann. Probab. {\bf 37} (2009), no. 1, 206--249.

\bibitem{OrsingherPolito2012} E. Orsingher and F. Polito, The space-fractional Poisson process, Statist. Probab. Lett. {\bf 82} (2012), no. 4, 852--858.

\bibitem{Orsingher2017} E. Orsingher and B. Toaldo, Space-time fractional equations and the related stable processes at random time,
J. Theoret. Probab. {\bf 30} (2017), no. 1, 1--26.

\bibitem{Pruss} J. Pr\"uss, Evolutionary Integral Equations and Applications, Birkh\"auser, Basel, 1993.

\bibitem{RGM2025} S. V. Rogosin, F. Giraldi and F. Mainardi, On differentiation with respect to parameters of the functions of the Mittag-Leffler type, Integral Transforms Spec. Funct. {\bf 36} (2025), no. 5, 400--412.

\bibitem{SKM1993} S. G. Samko, A. A. Kilbas and O. I. Marichev, Fractional Integrals and Derivatives. Theory and Applications, Gordon and Breach, Yverdon, 1993.

\bibitem{Schilling2012} R. L. Schilling, R. Song and Z. Vondracek, Bernstein Functions: Theory and Applications, 2nd ed., De Gruyter Studies in Mathematics {\bf 37}, De Gruyter, Berlin, 2012.

\bibitem{StanleyEC1} R. P. Stanley, Enumerative Combinatorics. Vol. 1, second edition, Cambridge University Press, Cambridge, 2012.

\bibitem{Stinga2019} P. R. Stinga, User's guide to the fractional Laplacian and the method of semigroups, in: Handbook of Fractional Calculus with Applications, vol. 2, De Gruyter, Berlin, 2019, pp. 235--265.

% \bibitem{Turban2004} L. Turban, Iterated random walk, Europhys. Lett. {\bf 65} (2004), no. 5, 627--632.

\bibitem{Ulgenes2025} D. P. H. Ulgenes, Series and product representations of gamma, pseudogamma and inverse gamma functions, Maple Transactions {\bf 5} (2025), no. 1, Article 16090.
\bibitem{VanMieghem2021} P. Van Mieghem, The Mittag-Leffler function, arXiv:2005.13330v4, 2021.
\end{thebibliography}
\end{document}